\def\IH{{\Bbb H}} 
 
\def\IR{{\Bbb R}} 
 
\def\IS{{\Bbb S}} 
 
\def\IK{{\Bbb K}}

\def\IC{\Bbb C} 
\def\ID{{\Bbb D}}

\def\mP{\mathcal P}
\def\mE{\mathcal E}
\def\mF{\mathcal F}
\def\mI{\mathcal I}

\def\oD{{\overline{\ID}}}

\def\zbar{{\overline{z}}} 
 
\def\xibar{{\overline{\xi}}} 
\def\wbar{{\overline{w}}} 
\def\mubar{{\overline{\mu}}}

\documentclass[12pt]{article} 

\usepackage[psamsfonts]{amssymb} 
\usepackage{amsfonts,amsmath,mathrsfs}  
\usepackage{epsfig,multicol}

\newtheorem{theorem}{Theorem}
\newtheorem{lemma}{Lemma}
\newtheorem{corollary}{Corollary}
\newtheorem{definition}{Definition}
\numberwithin{equation}{section}

\usepackage{chngcntr}
\usepackage{apptools}
\AtAppendix{\counterwithin{lemma}{section}}

\usepackage{hyperref}
\hypersetup{
    colorlinks,
    citecolor=black,
    filecolor=black,
    linkcolor=black,
    urlcolor=black
}

\title{The exponential Teichm\"uller theory: Ahlfors–Hopf differentials and diffeomorphisms}
\author{Gaven Martin \& Cong Yao \thanks{Work of both authors partially supported by the New Zealand Marsden Fund. C. Yao is also partially supported by the National Natural Science Foundation of China (No. 12401096), the Natural Science Foundation of Shandong Province (No.ZR2024QA035).
}
}
\date{}
\begin{document}
\maketitle

\begin{abstract}
We consider the minimisers of the $p$-exponential conformal energy for homeomorphisms $f:R \to S$ of finite distortion $\IK(z,f)$ between analytically finite Riemann surfaces in a fixed homotopy class of a barrier function $f_0$, 
\[ \mE_p(f:R,S)=\int_R \exp(p\IK(z,f))\;  d\sigma(z). \]
Homeomorphic minimisers always exist should the barrier function be a homeomorphism of finite energy,  $\mE_p(f_0,R,S)<\infty$.  Lifting to the universal covers leads us to the weighted $p$-exponential conformal energy   
\[ \int_\mP \exp(p\IK(z,f))\; \eta(z) \, dz, \]
where $\mP\subset\ID$ is a fundamental polyhedron, $\eta(z)$ is the hyperbolic metric,  along with an automorphy condition for $f$. This problem is {\em not variational}, however the Euler-Lagrange equations show that the inverses $h=f^{-1}$ of {\em sufficiently regular} stationary solutions have an associated holomorphic quadratic differential -- the Ahlfors-Hopf differential, 
\[ \Phi=\exp(p\IK(w,h))\,h_w\overline{h_\wbar}\,d\sigma_R(h). \]
When lifted to the disk the Ahlfors-Hopf differential has the form
\[ \Phi=\exp(p\IK(w,h))\, h_w\overline{h_\wbar}\,\eta(h).\] 
As a consequence of the Riemann-Roch theorem and an approximation technique,  we show that the variational equations do indeed hold for extremal mappings of finite distortion $f:R \to S$ between analytically finite Riemann surfaces. We then take this differential as a starting point for higher regularity and establish that if $h:\Omega\to\tilde{\Omega}$ is a Sobolev homeomorphism between planar domains with
a holomorphic Ahlfors-Hopf differential,    
then $h$ is a diffeomorphism.  It will follow that $h$ is harmonic in a metric induced by its own (smooth) distortion.  Key parts of the proof are developing equations for the Beltrami coefficient of a mapping with holomorphic Ahlfors-Hopf differential, and this is done by establishing a connection between associated degenerate elliptic non-linear Beltrami equations and these harmonic mappings.  Another key fact is an analogue of the higher regularity theory for solutions of distributional inner variational equations established in our earlier work. 

We then pull the problem back to the surfaces to conclude that minimisers $f_p\in [f_0]$ of $\mE_p(f:R,S)$ are diffeomorphisms and are unique stationary points. This now allows us to link two different approaches to Teichm\"uller theory,  namely the classical theory of extremal quasiconformal maps and the harmonic mapping theory.  As $p\to\infty$ we show $f_p\to f_\infty$ to recover the unique extremal quasiconformal mapping first rigourously proved by Ahlfors.  This extremal quasiconformal mapping is not a diffeomorphism (unless it is conformal) and $f_p$ degenerate on a divisor associated with the associated holomorphic quadratic differential of $f_\infty$.  As $p\to0$ we recover the harmonic diffeomorphism in  $[f_0]$ and Shoen-Yau's results.

In the last part we give some remarks on the the generalised problem for mappings of exponential finite distortion. In this case the existence of the Ahlfors-Hopf differential is implied if, for instance, $\exp(q\IK(zf))\in L^1_{loc}(\ID)$ for any $q>p$,  or if a Hamilton sequence condition is satisfied.
\end{abstract}

\newpage

\tableofcontents

\section{Introduction}
 
\subsection{Quasiconformal maps and Teichm\"uller Theory.} In 1939 Teichm\"uller stated his famous theorem \cite{T}: In the homotopy class of a diffeomorphism between closed Riemann surfaces, there is a unique extremal quasiconformal mapping of smallest maximal distortion 
\begin{equation}\label{Kdef}
\IK(z,f) = \frac{\|Df\|^2}{J(z,f)} = \frac{|f_z|^2+|f_\zbar|^2}{|f_z|^2-|f_\zbar|^2}.
\end{equation} 
Furthermore, this function has Beltrami coefficient 
\begin{equation}\label{mudef}
\mu_f=f_\zbar/f_z 
\end{equation}
 of the form
\begin{equation}\label{teich1}
\mu_f(z) =k\;  \frac{\overline{\Phi(z)}}{|\Phi(z)|},
\end{equation}
where $\Phi$ is a holomorphic function and $0\leq k <1$ is a constant.  The function $\Phi$ defines a holomorphic quadratic differential and gives local coordinates away from its zero set in which the extremal map is linear.

\medskip

This, together with Ahlfors' 1935 work \cite{Ahlfors1} established quasiconformal mappings in function theory.  Teichm\"uller's pioneering work contains a lot of conjectures and incomplete proofs. Later in 1953 Ahlfors reconsidered that work in  \cite{Ah3} and gave a complete proofs and a systematic introduction to this topic.  Ahlfors' approach was through the Calculus of Variations.  He posits the problem in the Universal Teichm\"uller space of quasiconformal mappings of the unit disk $\ID$,  constrained by an automorphy condition determined by the homotopy class of the surface diffeomorphism.   

\medskip

There are no variational equations that a minimiser of the maximal distortion satisfies,  so Ahlfors sets up the natural problem to identify the minimiser $h$ of the $L^p$-mean distortion. Then,  when viewed via the hodographic transformation (that is Ahlfors considers the related functional defined on the inverses), he computes the inner-variational Euler-Lagrange distributional equations to find  that a non-linear combination of the first order weak derivatives
\begin{equation}\label{AhlforsHopf} \Phi= \IK(w,h)^{p-1} h_w\overline{h_\wbar}, \end{equation}
is holomorphic (see \cite[\S 4, pg. 45]{Ah3}), and so we call the function $\Phi$ defined at (\ref{AhlforsHopf}) the Ahlfors-Hopf differential. Now Ahlfors analyses what happens as $p\to\infty$ to obtain the Teichm\"uller mapping (the function $\Psi$ in (\ref{teich1}) is the $L^1(\ID)$ normalised limit of the holomorphic Ahlfors-Hopf differentials - this limit is non-vanishing by the Riemann-Roch theorem.  

\medskip

In many ways Ahlfors work is prescient.  The technology of the times did not allow Ahlfors to say very much at all about the $L^p$ minimisers (or more precisely their inverses).  In fact he was working within the now well established framework of mappings of finite distortion \cite{IM1,AIM} and gave perhaps the first key modulus of continuity estimates for such mappings,  so as to establish continuity of a minimiser,  see \cite[Theorem 12]{Ahlfors1}.   It seems almost a complete miracle that the limit of these $L^p$-minimisers,  none of which we know to be in any way regular -- let alone homeomorphic,  is a quasiconformal mapping. 

\medskip

\subsection{Harmonic maps and Teichm\"uller Theory.} Next,  from around 1989 there has been another approach to Teichm\"uller theory through harmonic mappings -- initiated in the important work of Wolf \cite{Wolf} and Tromba   \cite{Tromba} but based on much earlier work,  perhaps most notable in our context is the theorem of Schoen-Yau \cite{SY} showing harmonic mappings between closed Riemann surfaces are diffeomorphisms.   We recommend \cite{DW, Jost} for the basic theory here.   In this setting of harmonic mappings there are equations from the outset -- the harmonic mapping equation,  also known as the {\em tension equation}.  There are natural associated geometric structures -- for instance the Weyl-Petersson metric on Teichm\"uller space,  see Wolpert's survey, \cite{Wolpert}.  Considerable advances over the last few decades have now shown that both approaches can be used to establish most of the important results in the area.  As a nice example is the Nielsen realisation problem which has proofs via both approaches, see Kerckhoff \cite{Kerckhoff}, Wolpert \cite{Wolpert2} and Tromba \cite{Tromba2}.

\medskip 

Each of these approaches -- harmonic and quasiconformal -- has its virtues and drawbacks.  For instance the Teichm\"uller metric is complete,  but only ``almost negatively curved'',  whereas the Weil-Petersson metric is negatively curved but seldom complete.  Harmonic mappings are diffeomorphisms but extremal quasiconformal mappings are not.

 \medskip

The approach we offer here seeks to unify both the harmonic and quasiconformal approaches as various limiting regimes ($p\to\infty$ and $p\to 0$) of the space of extremal mappings of exponential distortion,  that is minimisers of $\mE_p(\cdot)$, given the automorphy condition (or boundary data).  Already in \cite{AIMO,IMO} loose connections were established in the $L^p$ case and these were refined to establish existence and higher regularity in our work \cite{MY1,MY3}.  We need some definitions to clearly explain what's going on here.\\

\subsection{Mappings of finite distortion.}  The main objects of our study in this article are the finite distortion mappings.  The basic references here are the texts \cite{AIM,HK,IM1}.

\medskip

Let $\Omega\subset\IC $ be a planar domain and $f\in W_{loc}^{1,1}(\Omega)$, the Sobolev space of functions whose first derivatives are locally integrable functions. Let $Df$ denote the differential of $f$ and $J(z,f)=\det(Df)$ the Jacobian determinant.   Suppose that $J(z,f)\in L^{1}_{loc}(\ID)$ and that there is a measurable function $K(z)$ such that
\begin{equation}\label{1.5}
|Df(z)|^2\leq K(z) J(z,f), \hskip15pt \mbox{almost everywhere in $\Omega$.}
\end{equation}
Then, we say $f$ is a mapping finite distortion.  

Note the use of the operator norm in (\ref{1.5}). We define another distortion function  $\IK(z,f)$ by
\begin{equation}
\IK(z,f)=
\begin{cases}
\frac{\|Df(z)\|^2}{J(z,f)} &J(z,f)\neq0\cr
1 &\mbox{otherwise}
\end{cases}.
\end{equation}
The following lemma is a simple calculation.
\begin{lemma}
\[ \IK(z,f) = \frac{1}{2}\left( {\bf K}(z,f)+\frac{1}{ {\bf K}(z,f)} \right),\]
where $ {\bf K}(z,f)$ is the smallest function for which (\ref{1.5}) holds. 
\end{lemma}

Thus mappings of finite distortion satisfy  the distortion inequality
\begin{equation}\label{distineq}
\|Df(z)\|^2\leq   \IK(z,f) J(z,f), \hskip15pt \mbox{almost everywhere in $\Omega$.}
\end{equation}
The function $ {\bf K}(z,f)$ is the standard distortion function used in much of the literature on quasiconformal mappings - a homeomorphism of finite distortion with $\| {\bf K}(z,f)\|_{L^\infty(\Omega)}\leq K<\infty$ is called $K$-quasiconformal.  The distortion $ {\bf K}$  does not have good convexity properties.  As Ahlfors already noted in his proof of Teichm\"uller's theorem, the fact that there is not a tangent at the identity is a real impediment to any approach by the calculus of variations,  though he describes it as a ``technical device''.  That is why he turned to the distortion function $\IK(z,f)$ (see \cite[\S V (30)]{Ahlfors1}).  For example,  consider the family of mappings
$(x,y)\mapsto (\lambda x,y/\lambda),\hskip15pt \frac{1}{2} \leq \lambda\leq 2$.
Then, 
\[  {\bf K}_\lambda= \max\{\lambda,1/\lambda\}, \hskip80pt \IK_\lambda = \frac{1}{2}\left( \lambda+\frac{1}{\lambda}\right).\]

\begin{center}
\scalebox{0.6}{\includegraphics*[viewport=-70 550 780 750]{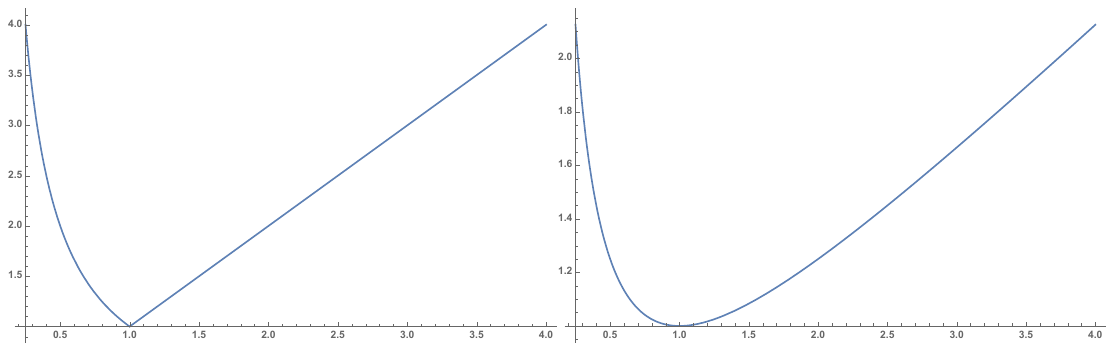}} \\
$ {\bf K}$ versus $\IK$ on a smooth family of diffeomorphisms.
\end{center}
 
In complex notation we have
\[ \|Df(z)\|^2 = |f_z|^2+|f_\zbar|^2, \hskip10pt
J(z,f) = |f_z|^2-|f_\zbar|^2.
\]
Thus 
\begin{equation} 
\IK(z,f)=\frac{1+|\mu_f(z)|^2}{1-|\mu_f(z)|^2}, \label{belt} \hskip10pt {\rm and}\hskip10pt
\mu_f(z) = \frac{f_\zbar(z)}{f_z(z)},
\end{equation}
is the Beltrami coefficient of $f$.

\medskip

It is not obvious,  but nevertheless a fact that if $f$ is a local homeomorphism in $W^{1,1}_{loc}(\Omega)$ the the Jacobian is locally in $L^1$, \cite{AIM,HK}.  Thus a $W^{1,1}_{loc}(\Omega)$ homeomorphism with $\IK(z,f)\in L^q_{loc}(\Omega)$ for some $q>0$ is a mapping of finite distortion.\\

A mapping has $p$-exponential finite distortion if 
\begin{equation}\label{energy}
\int_\Omega\exp(p\IK(z,f)) \; dz < \infty.
\end{equation}
There is a remarkably complete theory for these mappings which we now give a synopsis of.

\subsection{Mappings of exponentially integrable distortion.}\label{expsec}

The Orlicz-Sobolev space $W^{1,Q}_{loc}(\IC)$ consists of Sobolev functions whose weak first derivatives are $L^Q_{loc}(\Omega)$ functions.  In most of what follows
 \begin{equation} \label{kuu}
Q(t)=\frac{t^2}{\log(e+t)}.
\end{equation}
Thus our function spaces are just a little below the ``natural'' space $W^{1,2}_{loc}$ where local integrability of the Jacobian is assured.

\medskip
The following is the primary motivation for considering mappings with exponentially integrable distortion.  Suppose $f \in W^{1,1}_{loc}(\IC)$ is an orientation-preserving mapping
(that is,  $J(z,f)  \geq 0$) whose distortion function $\IK(z,f)$ satisfies $e^{p\IK(z,f)} \in L^1_{loc}(\IC)$ for some $ 0 < p < \infty$.  Then the inequality (\ref{1.5}) together with the elementary inequality
$
ab  \leqslant a \log(1+a) + e^b-1
$
 shows us
\[ \frac{\|Df\|^2}{\log(e+\|Df\|^2)}  \leq  \frac{\IK J}{\log(e+\IK J)} \leq   \frac{1}{p} \; \frac{J}{\log(e+J)} \; p \IK
\leq  \frac{1}{p}\left( J + e^{p\IK}-1\right)
\]
for all $p>0$.  Thus
\begin{equation} 
\int_\Omega \frac{\|Df\|^2}{\log(e+\|Df\|)} \leq \frac{2}{p}\int_\Omega J(z,f) + \frac{2}{p}\int_\Omega
[e^{p\IK(z,f)}-1]\; dz
\end{equation}
for any bounded domain.  These estimates directly imply that $f$ belongs to the Orlicz-Sobolev class $W^{1,Q}_{loc}(\IC)$ as soon as $J(z,f) \in L^1_{loc}(\Omega)$.

   Conversely, it is known  that for an orientation-preserving mapping $f$  the $L^2 \log^{-1}L$-integrability of the differential  implies that $J(z,f) \in L^1_{loc}(\Omega)$.   This slight gain in the regularity of the  Jacobian determinant is precisely why it will be possible to study solutions to the Beltrami equation using the space $W^{1,Q}_{loc}$.
 
 \medskip
 
One of the main tools in the theory is the generalization of the measurable Riemann mapping theorem and Stoilow factorization to the case of mappings with exponentially integrable distortion.

\medskip

\noindent \framebox[1.1\width]{\bf Existence of Solutions for Exponential Distortion}

\medskip

We slightly modify the statement of \cite[Theorem 20.4.9]{AIM} for our purposes.
\begin{theorem}\label{Existence} 
Let $\mu:\IC\to\ID$ be measurable.  Suppose there is $R>0$ such that $\mu(z)=0$ for $|z|\geq R$ and that the distortion function 
\begin{equation}
\IK(z) = \frac{1+|\mu|^2}{1-|\mu|^2}
\end{equation}
satisfies
\begin{equation}
\int_\IC e^{p\IK(z)}-e^p \; dz < \infty.
\end{equation}
Then the Beltrami equation 
\begin{equation} f_\zbar(z)=\mu(z) \; f_z(z)
\end{equation}
 admits a unique homeomorphic solution $f$ which is holomorphic near $\infty$ and has expansion
 \[ f(z) = z + a_1 z^{-1} + a_2 z^{-2}+ \cdots \]
 near $\infty$.  Further, 
\begin{equation} \label{varsQ}
f\in W^{1,Q}_{loc}(\IC), \quad \quad Q(t)=t^2 \log^{-1}(e+t).
\end{equation} 
Moreover,  every other $W^{1,Q}_{loc}$-solution $g$ to this  Beltrami equation in a domain $\Omega \subset \IC$ admits the factorization $$g = \phi \circ f,$$ where $\phi$ is a holomorphic function in the domain  $f(\Omega)$.
\end{theorem}

\medskip
\noindent \framebox[1.1\width]{\bf Optimal regularity.}

\medskip

Building on earlier work,  see \cite{IKM},  Astala, Gill, Rohde and Saksman established the optimal regularity for planar mappings with exponentially integrable distortion in \cite{AGRS}.

\begin{theorem} Let   $\Omega\subset\IC$ be a planar domain and suppose that the distortion function $\IK(z,f)$ of a mapping of finite distortion $f:\Omega\to \IC$ has $e^{p\IK(z,f )}\in L^{1}_{loc}(\Omega)$ for some $p > 0$. Then we have for every $0 < \beta < p$,
\begin{enumerate}
\item $\|Df \|^2 \log^{\beta-1}( e + \|Df \|) \in  L^{1}_{loc}(\Omega) $;
\item $J(z,f)\log^\beta(e+J(z,f)) \in  L^{1}_{loc}(\Omega)$.
\end{enumerate} 
Moreover, for every $p > 0$ there are examples that satisfy the hypotheses and yet fail both conclusions for $\beta = p$.
\end{theorem}

\medskip

\noindent \framebox[1.1\width]{\bf  Modulus of continuity.}

\medskip

Mappings of finite distortion with exponentially integrable distortion have continuous representatives.  In fact from \cite{GV, IM1}, and in particular \cite[Theorem 2.3 \& 2.4]{HK} we have the following theorem.

\begin{theorem}\label{thm4} Let $\Omega$ be a planar domain and $f\in W^{1,2}_{loc}(\Omega)$ be a mapping of finite distortion.  Then $f$ has a continuous representative.  

The same is true if $f\in W^{1,1}_{loc}(\Omega)$ has finite distortion and there is $p>0$ and $e^{p\IK(z,f)}\in L^{1}_{loc}(\Omega)$.
\end{theorem}

The best known modulus of continuity estimate is that of Onninen and Zhong \cite{OZ}.  Earlier results were established by David \cite{David}, so also \cite{IM1} and \cite{KO}.

\begin{theorem} \label{cont1} Let $f : \Omega \to \IC$ be a mapping of finite distortion whose distortion function $\IK(z,f)$ satisfies, for some $p > 0$,  
\[ K:=  \int_B e^{p\IK(z,f)} \; dz<\infty ,\]
where $B=\overline{\ID(z_0,R)}$  is a compact subset of $\Omega$.  Then, for every 
\[ z,w\in \ID\Big(z_0, \frac{R^e}{240^e}\Big[\frac{K}{\pi}\Big]^{(1-e)/2}\Big)
\]
we have
\begin{equation} |f(z)-f(w)| \leq \frac{C}{\log^{p/2}\big[\frac{K}{\pi|z-w|^2}\big]} \Big[\int_B J(z,f) \; dz\Big]^{1/2}.
\end{equation}
The constant $C$ depends on $R, K$ and $p$.
\end{theorem}

Following on from these we consider differentiability.  Open mappings in the Sobolev space $W^{1,1}_{loc}(\Omega)$ are differentiable almost everywhere by the Gehring-Lehto Theorem, \cite{GL}.  In particular homeomorphisms of finite distortion are differentiable almost everywhere.  However,  again more is true,  \cite[Corollary 2.25]{HK}.

\begin{theorem}  Let $p>0$ and let $f:\Omega\to\IC$ be a mapping of finite distortion with $e^{p \IK(z,f)}\in L^{1}_{loc}(\Omega)$.  Then $f$ is differentiable almost everywhere.
\end{theorem}

\noindent \framebox[1.1\width]{\bf Properties of inverses.}

\medskip

As noted in our discussion of Ahlfors' earlier work,  the properties of inverses,  and more precisely the equations they satisfy will be important.
We recall Theorem 1 in of Gill \cite{Gill},  building on earlier work of Hencl and  Koskela \cite{HK}.
\begin{theorem}\label{Gill} Let $\Omega$ be a planar domain and suppose $f:\Omega\to \IC $  is a homeomorphism of finite distortion and suppose also that $e^{p\IK(z,f)}\in L^1_{loc}(\Omega)$  for some $p > 0$. Then the inverse $h= f^{-1}:f(\Omega)\to \Omega$ is also a homeomorphism of finite distortion with associated distortion function $\IK(w,h)$ and
\begin{equation}
\IK(w,h)\in L^\beta_{loc}(\Omega), \hskip15pt \mbox{for all $0 < \beta < p$.}
\end{equation}
This result is sharp in the sense that there exist homeomorphisms $f$ as above for which $\IK(w,h)\not\in L^p_{loc}(f(\Omega))$.  
\end{theorem}

Notice that if $f:\Omega\to f(\Omega)=\Omega' $  is a homeomorphism of finite distortion with $e^{p\IK(z,f)}\in L^1_{loc}(\Omega)$  for some $p > 0$ and if $\IK(w,h)\in L^1(\Omega')$,  then (using some of the following results - eg condition ${\cal N}$ and ${\cal N}^{-1}$ and almost everywhere differentiability) we compute that
\[ \int_{\Omega'} \IK(w,h) \; dw = \int_{ \Omega} \IK(f,h) \; J(z,f) dz = \int_\Omega\|Df(z)\|^2 \; dz  \]
so that in fact $f\in W^{1,2}(\Omega)$.  This easy example shows how mean distortion on one hand and Dirichlet energy on the other are connected,  see \cite{AIM,HK}.

\bigskip
 
\noindent \framebox[1.1\width]{\bf  Boundary values and extension.} 

\medskip

With a little work the results above quickly yield the following theorem which appears a little disguised in work of Ryazanov, Srebro and Yakubov \cite{RSY}.  Within our framework \cite{KVO} proves a little more.

\begin{theorem}\label{extend}  Let $f:\ID\to\ID$ be a homeomorphism (onto) and suppose that $e^{p\IK(z,f)}\in L^1(\ID)$ for some $p>0$.  Then $f$ extends continuously to  homeomorphism $f:\ID\to\ID$ and then by reflection to a mapping of finite distortion $\tilde{f}:\IC\to\IC$ with $e^{p\IK(z,\tilde{f})}\in L^{1}_{loc}(\IC)$.
\end{theorem}

This leads us to the following cleaner modulus of continuity estimate which we will use.

\begin{corollary}\label{Extension}
Let $f:\ID\to\ID$ be a homeomorphism (onto) and suppose that $e^{p\IK(z,f)}\in L^1(\ID)$ for some $p>0$.  Then there is a constant $C_p$,  depending only on $p$, $|f(0)|$ and $\int_\ID e^{p\IK(z,f)}\;dz$,  so that we have the modulus of continuity estimate :  for all $z,w\in\ID$, 
\begin{equation} |f(z)-f(w)| \leq \frac{C_p}{\log^{p/2} \frac{1}{|z-w|} } .
\end{equation}
\end{corollary}
\noindent{\bf Proof.}  Let $z_0=f(z)$ and put $\phi(z)=(z-z_0)/(1-\overline{z_0} z)$.    Then as a M\"obius transformation of the disk we have $\IK(z,\phi\circ f)=\IK(z,f)$.  Also $\phi$ is a bilipschitz homeomorphism of $\overline{\ID}$ whose bilipschitz constant depends on $|f(0|$ and is finite:
\begin{eqnarray*}
|\phi(z)-\phi(w)| &= & \left| \frac{z-z_0}{1-\overline{z_0} z}-\frac{w-z_0}{1-\overline{z_0} w} \right| \leq   \frac{1+|z_0|}{1-|z_0|} \; |z-w|.
\end{eqnarray*} 
A similar lower bound can be found as well.  We can now absorb this Lipschitz constant into $C_p$ and assume that $f(0)=0$.
Then we observe that first $f$ extends to a homeomorphic mapping of the closed disk. Then we can extend by reflection to a the complex plane by
\[ 
\tilde{f} (z):=\begin{cases}
f (z)&z\in\overline{\ID }\\
\frac{1}{\overline{f (\frac{1}\zbar )}}&z\in \IC \setminus \overline{\ID }
\end{cases}.
\] 
Let $R>1$ and put $r=1/R$.  We can compute, for $z\in\ID(0,R)\setminus \overline{\ID }$,
\[
(\tilde{f})_z(z)=\frac{1}{\overline{f^2(\frac{1}\zbar )}}\overline{(f)_z(\frac{1}\zbar )\frac{1}{\bar{z}^2}}, \;\; 
(\tilde{f})_\zbar (z)=\frac{1}{\overline{f^2(\frac{1}\zbar )}}\overline{(f)_\zbar (\frac{1}\zbar )\frac{1}{z^2}}, 
\]
and so $ \IK(z,\tilde{f})=\IK(\frac{1}\zbar ,f)$.  
In particular, $\tilde{f}_0$ has finite distortion in $A=\ID(0,R)\setminus\overline{\ID}$, and for any continuous $\Psi$
\begin{eqnarray*}
\int_{A}\Psi(\IK(z,\tilde{f}_0))\; dz&=&\int_{A}\Psi(\IK(\frac{1}\zbar ,f_0)\; dz =\int_{\ID -\ID(0,r)}\Psi(\IK(\zeta,f_0))\frac{1}{|\zeta|^4}\; d\zeta\\
&\leq&\frac{1}{r^4}\int_\ID \Psi(\IK(\zeta,f_0))\; d\zeta.
\end{eqnarray*}
With $\Psi(t) = e^{pt}$ we achieve 
\[ \int_{\ID(0,R)} e^{p\IK(z,\tilde{f})} \; dz \leq 1+R^4\int_{\ID} e^{p\IK(z,f)}\; dz. \] 
The result follows directly from Theorem \ref{cont1} as soon as we choose $R$ large enough so that $\frac{R^e}{240^e}\Big[\frac{K}{\pi}\Big]^{(1-e)/2}\geq 2$. \hfill $\Box$

\bigskip
A moments thought will convince that these modulus of continuity estimates must depend on $|f(0)|$ as they must do so for conformal mappings.  

\bigskip

\noindent \framebox[1.1\width]{\bf The area formulae and Lusin's condition $\mathcal{N}$.}

\medskip

The area formulae deal with the question of change of variables.  For Sobolev mappings we will need the following Lusin condition $\mathcal{N}$:

\begin{definition}
Let $f:\Omega\to\IC$ be measurable. We say $f$ satisfies Lusin's condition $\mathcal{N}$, if for any $E\subset\Omega$ with $|E|=0$, we have $|f(E)|=0$. If $f$ is a homeomorphism. We say $f$ satisfies Lusin's condition $\mathcal{N}^{-1}$, if its inverse $f^{-1}$ satisfies Lusin's condition $\mathcal{N}$.
\end{definition}
The next theorem follows directly from work of Gol'dstein and Vodopyanov \cite{GV} and Kauhanen,  Koskela and Mal\'y \cite{KKM1}.
\begin{theorem}
Let $p>0$ and let $f :\Omega\to \IC$ be a homeomorphism of finite distortion with $e^{p\IK(z,f)}\in L^{1}_{loc}(\Omega)$.  Then $f$ satisfies both $\mathcal{N}$ and  $\mathcal{N}^{-1}$.  In particular for all $\eta\in L^{1}(\Omega)$,
\[
\int_\Omega\eta(f(x)) \; J(x,f)dx=\int_{f(\Omega)}\eta(y) \; dy. 
\]
\end{theorem}

\noindent \framebox[1.1\width]{\bf Compactness.}

\medskip
 
With the above results at hand it is not too hard to see how the following refinement of \cite[Theorem 8.14]{IM1} is proved.  They only key ideas missing from the discussion above is the convexity of the function $e^{p\IK(z,f)}$ as a function of the minors of $Df$ (connecting with Ball's notion of polyconvexity) and  a study of  distributional Jacobians and their weak convergence.   The following result will easily provide us with the existence of minimisers.
 
 \begin{theorem}\label{compact}  Let $\Omega\subset\IC$ be a planar domain.  Given $p>0$, $M<\infty$ and $z_0,z_1\in \IC$ let ${\cal F}$ denote the family of mappings $f:\Omega\to\IC$ of finite distortion such that
 \begin{itemize}
 \item $f(\Omega)\subset \IC\setminus \{z_0,z_1\}$,
 \item \[ \int_\Omega e^{p\IK(z,f)} \; dz \leq M. \]
 \end{itemize}
 Let $U\subset \Omega$ be relatively compact.  Then
 \begin{enumerate}
 \item ${\cal F}$ is bounded in $W^{1,Q}(U)$.  $Q(t)=t^2\log^{-1}(e+t)$.
 \item ${\cal F}$ is closed under weak convergence in $W^{1,Q}(U)$.
 \item ${\cal F}$ is equicontinuous in $U$.
 \item The limit of a locally uniformly convergent sequence of mappings in ${\cal F}$ belongs to ${\cal F}$.
 \item The limit of a locally uniformly convergent sequence of homeomorphisms in ${\cal F}$ is a homeomorphism in ${\cal F}$ if it is non-constant.
 \end{enumerate}
 \end{theorem}
 \noindent {\bf Proof.}  There is a conformal mapping $\Phi:\IC\setminus \{z_0,z_1\}$ and for any $f\in{\cal F}$ we obviously have
 \[ \int_\Omega e^{p\IK(z,\Phi\circ f)} \; dz=\int_\Omega e^{p\IK(z, f)} \; dz \leq M, \]
 and so we could assume the family ${\cal F}$ is bounded and the first four conclusions of  Theorem \ref{compact} follow from \cite[Theorem 8.14]{IM1}.  Suppose $\{f_j\}$ is locally uniformly convergent sequence of homeomorphisms in ${\cal F}$.  Then by 3. we may assume that this sequence has a non-constant limit $f\in {\cal F}$ and that $f\in W^{1,Q}(U)$. Let $V\subset \Omega$ be relatively compact and let $\mu_f$ be the Beltrami coefficient of $f$.  Define $\nu(z)$ as follows.  If $z\in V$ let $\nu(z)=\mu_f(z)$.  Otherwise put $\nu(z)=0$.  Since $f$ is continuous,  $f(\bar U)$ is compact in $f(\Omega)\subset \IC$.  Thus $\nu$ satisfies the hypotheses of Theorem \ref{Existence}.  Let $F$ be the unique homeomorphic (principal) solution described in that result.  Theorem \ref{Existence} now asserts that there is a conformal mapping $\varphi$ defined on $F(U)$ so that $f(z) = \varphi\circ F)(z)$.  However, as a uniform limit of homeomorphisms $f$ has topological degree equal to $1$ and is also a {\it monotone} mapping -- the preimage of a point is a connected set. Since $F$ is a homeomorphism we find that $\varphi$ is a conformal mapping of degree $1$ and so a local homeomorphism.  Thus $f$ is a local homeomorphism,  and monotone local homeomorphisms are  homeomorphisms. \hfill $\Box$

\subsection{Riemann surfaces and universal covers.}
Let $p>0$ and let $R$ be an analytically finite hyperbolic Riemann surface of signature $(g,n)$ - genus $g$ with $n$ punctures.  We realise $R=\IH/\Gamma$ where $\Gamma$,  isomorphic to the fundamental group of $R$, is a Fuchsian group of isometries of the hyperbolic plane $\IH$,  which we henceforth identify as the unit disk $\ID$ with hyperbolic metric $ds=(1-|z|^2)^{-1}|dz|$. Given a mapping $f$ between a pair of such surfaces,  say $R$ and $S$ we lift $f$ to the universal cover to obtain $\tilde{f}:\ID\to\ID$.  The covering projections $\pi_R:\ID\to R$ and $\pi_S:\ID\to S$ are locally conformal so we can define the point-wise distortion of $f$ as
\begin{equation}\label{1.19}
\IK(z,f)=\IK(w,\tilde{f}),  \hskip15pt w\in\ID, \;\; \pi_R(w)=z\in R,
\end{equation}
since 
\begin{equation}\label{1.20}
(\pi_S\circ\tilde{f})(w)=(f\circ \pi_R)(w).
\end{equation}
Then the definition of mappings of finite distortion is straightforward.

For any function $F:R \to \IC$ the integral with respect to the hyperbolic surface measure $\sigma(z)$ on $R$ is
\begin{equation}\label{1.21}
\int_R F(z) \; d\sigma(z) = \int_{\mP} F(\pi_R(w)) \; \frac{dw}{(1-|w|^2)^2},
\end{equation}
where $\mP\subset \ID$ is any convex fundamental polyhedron for the action of $\Gamma$ on $\ID$ -- there is always one such.

\medskip

We now define the conformal energies we consider for maps between Riemann surfaces.  Let $\Psi:[1,\infty)\to[1,\infty)$ be convex and increasing,  and $f:R\to S$.  Define
\begin{equation}\label{1.22}
{\cal E}_\Psi(f) = \int_R \Psi(\IK(z,f))\; d\sigma(z).
\end{equation}
The following theorem connects these energies with our earlier studies.
\begin{theorem} Let $f:R \to S$ be a mapping of finite distortion. Let $\tilde{f}:\ID\to\ID$ be a lift of $f$ to the universal covers $\ID$.  Then
\begin{equation}
{\cal E}_\Psi(f)  = \int_\mP \Psi(\IK(w,\tilde{f}))\; \frac{dw}{(1-|w|^2)^2} .
\end{equation}
which is independent of the choice of fundamental polyhedron $\mP$.
\end{theorem}
\noindent{\bf Proof.} The expression follows from (\ref{1.21}). To see it does not depend on the choice of $\mP$, we make the following calculation using the identity 
\[ \frac{|\gamma'(z)|}{1-|\gamma(z)|^2}=\frac{1}{1-|z|^2}\]
for a M\"obius transformation of the disk.  Then
\begin{eqnarray*}
\int_{\gamma(\mP)}\Psi(\IK(z,\tilde{f}))\;\frac{dz}{(1-|z|^2)^2}&=&\int_\mP\Psi(\IK(\gamma(w),\tilde{f}))\;\frac{|\gamma'(w)|^2dw}{(1-|\gamma(w)|^2)^2}\\
&=&\int_\mP\Psi(\IK(w,\tilde{\gamma}\circ\tilde{f}))\;\frac{dw}{(1-|w|^2)^2}\\
&=&\int_\mP\Psi(\IK(w,\tilde{f}))\;\frac{dw}{(1-|w|^2)^2}.
\end{eqnarray*}
This establishes the result.\hfill $\Box$

\medskip

We can now consider the problem of identifying the minimisers of mean distortion for mappings between Riemann surfaces (more precisely the functional ${\cal E}_\Psi$).  We simply restrict the class of maps we admit as candidates to those which are automorphic with respect to the Fuchsian groups $\Gamma$ and $\tilde{\Gamma}$ of $R$ and $S$ respectively.  Hence given a pair of analytically finite surfaces of the same signature, associated Fuchsian groups and an isomorphism $\mI:\Gamma\to \tilde{\Gamma}$ (all of which we fix henceforth,  but write $\mI(\gamma)=\tilde{\gamma}$ for brevity) we set
\begin{equation}\label{fproblem}
{\cal F}_\Psi = \big\{f:\ID\to\ID: f \circ \gamma = \tilde{\gamma}\circ f \;\; \& \;\; \int_\mP \Psi(\IK(z,f)) \; \frac{dz}{(1-|z|^2)^2}< \infty \big\}.
\end{equation}
It is the dual family we will pay particular attention to for the proof of our main results.  Set
\begin{equation}\label{hproblem}
{\cal H}_\Psi = \big\{h:\ID\to\ID: h \circ \tilde{\gamma} = \gamma \circ h \;\; \& \;\; \int_\mP \Psi(\IK(w,h))J(w,h)\; \frac{dw}{(1-|h|^2)^2} < \infty \big\}.
\end{equation}

\begin{lemma} The families ${\cal F}_\Psi$ and ${\cal H}_\Psi$ are non-empty.  Every element of ${\cal H}_\Psi$ is continuous and extends continuously to the boundary $\IS$ and all elements share the same quasisymmetric boundary values. 
\end{lemma}
\noindent{\bf Proof.} Since $R$ and $S$ are analytically finite with the same signature,  there is a $K$-quasiconformal mapping between them.  The lift of this map to $\ID$ is quasiconformal and automorphic and easily seen to have finite energy.   Since $\Psi$ is convex and increasing we have $\Psi(t)\geq c_0 t$ and so every element of ${\cal H}$  is a $W^{1,2}(\ID)$ mapping of finite distortion and so has a continuous representative.  The boundary values of every $h\in {\cal H}$ are determined uniquely by the isomorphism defining the automorphy condition -- that is the homotopy class of the induced maps between the surfaces. Since there is a quasiconformal mapping in ${\cal H}$ all mappings share these quasisymmetric boundary values. This proves the lemma. \hfill $\Box$

\medskip

Note that elements of ${\cal F}$ need not be continuous for general $\Psi$, however they too share boundary values by virtue of the automorphic property.

\medskip

The next thing to note is that the automorphy property is preserved under local uniform limits.  Thus the following is immediate from our earlier results.

\begin{theorem}\label{thm14} Let $p>0$ and $\Psi(t) \geq c_0\, e^{p t}$ for some positive constant $c_0$.  Then there is a homeomorphic minimiser $f:R \to S$ of the functional ${\cal E}_\Psi(\cdot)$.
\end{theorem}

\section{Exponential minimisers}
Here we state our main problem.\\

Let $R,S$ be analytically finite Riemann surfaces, $0<p<\infty$. Let $f_0:R\to S$ be a homeomorphism with finite exponential distortion $\mE_p(f_0)<\infty$, where
\[
\mE_p(f):=\int_R\exp(p\IK(z,f)) \; d\sigma_R(z).
\]
Let $[f_0]$ denote the homotopy class of mappings between $R$ and $S$ defined by $f_0$ and set 
\[ 
\mF_p:=\{f\in W_{loc}^{1,1}(R):\mE_p(f)<\infty  \mbox{ and $f\in [f_0]$ is a homeomorphism} \}.\]
We wish to find the minimisers of $\mE_p(f)$ in the class $\mF_p$.\\

Lifting to the unit disk $\ID$, we may state an equivalent problem as follows:\\

Let $R=\ID/\Gamma$, $S=\ID/\tilde{\Gamma}$, where $\Gamma$ and $\tilde{\Gamma}$ are Fuchsian groups, $\mP$, $\tilde{\mP}$ be convex fundamental polyhedra. Let $f_0:\ID\to\ID$ be a homeomorphism which  is automorphic with respect to $\Gamma$ and $\tilde{\Gamma}$, that is  there is a isomorphism $\mI:\Gamma\to\tilde{\Gamma}$ such that $\mI(\gamma)\circ f_0=f_0\circ\gamma$, for all $\gamma\in\Gamma$, and  $f_0$ has finite exponential distortion $\tilde{\mE_p}(f_0)<\infty$, where
\[
\tilde{\mE}_p(f):=\int_\mP\exp(p\IK(z,f)) \;\frac{dz}{(1-|z|^2)^2}.
\]
We define the class $\tilde{\mF}_p$, where we will seek a minimiser, to be those functions  $f:\ID\to \ID$
\begin{enumerate}
\item $f\in W_{loc}^{1,1}(\ID)$ is a homeomorphism. 
\item $\tilde{\mE}_p(f)<\infty$.  
\item $\mI(\gamma)\circ f=f\circ\gamma$ for all $\gamma\in\Gamma$.
\end{enumerate}
Then our problem becomes to find the minimisers of $\tilde{\mE}_p(f)$ in the class of $\tilde{\mF}_p$. As we have noted,  the automorphy condition implies $f$ extends to a homeomorphism $\overline{\ID}\to \overline{\ID}$ which quasisymmetric on the boundary.

From now on, since we have already established the equivalence between these two models, we  abuse notation  and do not distinguish between them.

\subsection{Variational equations.}
In the previous section we explained the existence of a minimiser in a given homotopy class of homeomorphisms. Now we seek equations satisfied by a minimiser (or stationary point) using the calculus of variations. We consider the disk model
\[
\mE_p(f)=\int_\mP\exp(p\IK(z,f)) \;\eta(z)dz,
\]
where $\eta(z)=\frac{1}{(1-|z|^2)^2}$. Let $\varphi\in C_0^\infty(\mP)$, $g^t=z+t\varphi$, $f^t=f\circ(g^t)^{-1}$. 

We say $f$ is {\it variational}, if $\mE_p(f^t)<\infty$ for a small variation $t\in(-\epsilon,\epsilon)$. A variational minimiser $f$ satisfies
\[
0=\frac{d}{dt}\Big|_{t=0}\int_\mP\exp(p\IK(z,f^t)) \; \eta(z)dz.
\]
Direct computation (following e.g. \cite[\S 21.2]{AIM} in the $L^p$ case) leads to the {\it inner variational equation }
\begin{equation}\label{2.1}
\int_\mP\exp(p\IK(z,f))(\eta\varphi)_z\; dz=\int_\mP\exp(p\IK(z,f))\frac{2p\overline{\mu_f}}{1-|\mu_f|^2}\eta\varphi_\zbar\; dz,
\end{equation}
for every $\varphi\in C_0^\infty(\mP)$. The same process applies to  the inverse function $h$ as well and this gives
\begin{equation}\label{2.2}
\int_{\tilde{\mP}}\exp(p\IK(w,h))h_w\overline{h_\wbar}\eta(h)\psi_\wbar\; dw=0,
\end{equation}
for every $\psi\in C_0^\infty(\tilde{\mP})$. Now Weyl's lemma gives that
\begin{equation}\label{2.3}
\Phi:=\exp(p\IK(w,h))h_w\overline{h_\wbar}\eta(h)
\end{equation}
is holomorphic in $\tilde{\mP}$, provided that $\Phi\in L^1(\tilde{\mP})$. Since this applies on an arbitrary convex fundamental polyhedron $\tilde{\mP}$ with respect to $\tilde{\Gamma}$, and since any point of $\ID$ lies in such a polyhedron, unique analytic continuation implies $\hat{\Phi}$ is a holomorphic function in $\ID$ and using the automorphy condition, it must satisfy the following equation:
\begin{equation}\label{2.4}
\Phi(w)=\Phi(\tilde{\gamma})(\tilde{\gamma}')^2,\quad\forall\tilde{\gamma}\in\tilde{\Gamma}.
\end{equation}
This also pulls back to $S$ as a holomorphic quadratic differential
\begin{equation}\label{2.5}
\Phi=\exp(p\IK(w,h))h_w\overline{h_\wbar}\;d\sigma_R(h).
\end{equation}
Note here $\|\Phi\|_{L^1(S)}=\|\Phi\|_{L^1(\tilde{\mP})}$, which is independent of the choice of $\tilde{\mP}$. 

As per \S1, this holomorphic function $\Phi$ is called the Ahlfors-Hopf differential and will play an important role in our study. In fact we will prove that $h$ is diffeomorphic based on $\Phi$ and (\ref{2.5}) in \S 2 \& 3. Another consequence of the existence of $\Phi$ will be to prove the uniqueness of minimisers, see \cite{MY5}.

We also remark that the variationality of stationary points is not automatically true, \cite{MY2}. An easy observation is that if $\mE_q(f)<\infty$ for  some $q>p$, then $f$ is variational. However, in the next section we will prove that $\Phi$ always exists in our case of maps between analytically finite surfaces.  We do not know what happens in more generality but speculate that it may not be the case that there is always an Ahlfors-Hopf differential.  Note that there is always such a differential in the case of $L^p$-minimisers,  but these are not always homeomorphisms \cite{MY1}. This presents a very interesting dichotomy. \\

\subsection{Existence of Ahlfors-Hopf differentials.}
As we mentioned above, the exponential minimisers are not automatically variational. However, we will prove there is a minimiser which has Ahlfors-Hopf differential anyway, and what we want to do is show this minimiser is a smooth diffeomorphism and that it is unique. The additional tool we have here is the Riemann-Roch theorem,  see Serre and Borel \cite{Serre}, but for a nice exposition see the notes by Kapovich \cite{Kap}. In fact it follows from the Riemann-Roch theorem that a closed Riemann surface of genus $g> 1$ has $3g -3$ linearly independent
holomorphic quadratic differentials over the field of complex numbers. Each
differential has $4 g - 4$ zeros. A similar result is true for analytically finite Riemann surfaces.\\

The $L^p$ problems are variational, and so for each $p$ there exits a holomorphic Ahlfors-Hopf differential, see \cite{IMO},\cite{MY1}. 

We therefore consider  the truncated problems
\[
\mathsf{E}_N(f):=\sum_{n=0}^N\int_R\frac{p^n\IK(z,f)^n}{n!}\;d\sigma_R=\sum_{n=0}^N\int_S\frac{p^n\IK(w,h)^n}{n!}\;d\sigma_R(h).
\]
As linear combinations of $L^p$ problems, following the arguments of \cite[\S 5 \& Theorem 5.2]{MY1}, we obtain for each $N$  a minimiser $f_N$ and a holomorphic quadratic differential
\[
\Phi_N=\sum_{n=0}^{N-1}\frac{p^n\IK(w,h_N)^n}{n!}(h_N)_w\overline{(h_N)_\wbar}\; d\sigma_R(h_N).
\]
where $h_N=f_N^{-1}$ is the pseudo-inverse minimiser of $\mathsf{E}_N$ in the enlarged space, see \cite[\S5 \& Theorem 5.1]{MY1}. By the Riemann-Roch theorem, the space of quadratic differentials is finite dimensional, so we get a non-degenerate sequence 

\begin{equation}\label{HS}
\frac{\Phi_N}{\|\Phi_N\|_{L^1(S)}}\not\to0,
\end{equation} and  thus $\Phi_N\to\Phi$ where $\Phi$ is also a non-zero holomorphic quadratic differential. We remark that such a sequence as at (\ref{HS}) is called a Hamilton sequence, \cite{EL}.  As with extremal quasiconformal mappings and boundary value problems,  the existence such a Hamilton sequence seems a crucial step in the boundary value problem for extremal mappings of exponential distortion to obtain smoothness and uniqueness.

Returning to our argument,  on  the other hand, the functions $f_N$, $h_N$ converge weakly in $W^{1,1}(R)$ and $W^{1,2}(S)$ respectively, to some $f$ and $h$. As we will see below, this $f$ actually has finite exponential distortion, so $f$ is a homeomorphism and $f=h^{-1}$.
\begin{lemma}\label{EConv}
Let $f_N$ be a minimiser of $\mathsf{E}_N$, $f_N\to f$ weakly in $W^{1,1}(S)$, then
\[
\mE_p(f)=\lim_{N\to\infty}\mathsf{E}_N(f_N).
\]
Moreover, $f$ is a minimiser in $\mF_p$.
\end{lemma}
\noindent{\bf Proof.} Note $\mathsf{E}_N(f_N)$ is a bounded non-decreasing sequence. Also, by Fatou's lemma and the polyconvexity of the functional we have
\[
\mE_p(f)\leq\lim_{N\to\infty}\mathsf{E}_N(f)\leq\lim_{N\to\infty}\mathsf{E}_N(f_N).
\]
However, if the inequality holds, then for sufficiently large $N$,
\[
\mathsf{E}_N(f)\leq\mE_p(f)<\mathsf{E}_N(f_N),
\]
which is impossible as $f_N$ is a minimiser of $\mathsf{E}_N$. The same contradiction occurs if the inequality holds for any other function $g\in\mF_p$, which proves the second claim.\hfill$\Box$

\bigskip

Lemma \ref{EConv} gives the Radon-Riesz property, and this implies that we have strong convergence of $f_N\to f$ and $h_N\to h$, see \cite{MY3}. We then obtain the holomorphic quadratic differential (\ref{2.5}). Furthermore, when $h$ is lifted to the disk, (\ref{2.3}) and (\ref{2.4}) hold. We record this as follows:
\begin{theorem}\label{AHdiff}
In the class $\mF_p$, there exists a homeomorphic minimiser $f$ whose inverse $h$ has holomorphic Ahlfor-Hopf differential (\ref{2.5}), and its lift satisfies (\ref{2.3}), (\ref{2.4}) in $\ID$.
\end{theorem}

\subsection{Diffeomorphic minimisers.}
In this subsection we consider the ``$f$-side''  inner variational equation (\ref{2.1}). Let 
\[
\sigma=-\mathcal{C}^*|_\ID\left(\exp\left(p\IK(z,f)\right)\eta_z\right).
\]
where $\mathcal{C}^*|_\ID$ is the adjoint Cauchy transform on $\ID$, see \cite[Section 4.1]{AIM}. This enable us to rewrite (\ref{2.1}) as
\begin{equation}\label{2.6}
\int_\mP\exp(p\IK(z,f))\frac{2p\overline{\mu_f}}{1-|\mu_f|^2}\eta\varphi_\zbar\;dz=\int_\mP\left((\exp(p\IK(z,f))-e^p)\eta+\sigma\right)\varphi_z\;dz,
\end{equation}
for every $\varphi\in C_0^\infty(\ID)$. If we assume $\IK(z,f)\exp(p\IK(z,f)) \in L_{loc}^2(\mP)$, then there exists an $F\in W_{loc}^{1,2}(\mP)$ such that
\begin{equation}\label{2.7}
F_z=\exp(p\IK(z,f))\frac{2p\overline{\mu_f}}{1-|\mu_f|^2}\eta,\quad F_\zbar=(\exp(p\IK(z,f))-e^p)\eta+\sigma.
\end{equation}
Thus
\begin{equation}\label{2.8}
F_\zbar=\mathcal{A}_p\left(\frac{|F_z|}{\eta}\right)\eta+\sigma,
\end{equation}
where $\mathcal{A}_p$ is the inverse of the function
\begin{equation}\label{2.9}
a_p(s)=(s+e^p)\sqrt{\log^2(s+e^p)-p^2},\quad s\geq0,
\end{equation}
We differentiate this to obtain
\[
a_p'(s)=\sqrt{\log^2(s+e^p)-p^2}+\frac{\log(s+e^p)}{\sqrt{\log^2(s+e^p)-p^2}}.
\]
Let $m_p=\min_{s\geq e^p}a_p'(s)$. It is not hard to see that $m_p>1$ for all $p>0$, and $m_p\to1$ as $p\to0$; $m_p=\mathcal{O}(\sqrt{p})$ as $p\to\infty$. Examples are given in the graphs below.
\begin{center}
\scalebox{0.4}{\includegraphics*{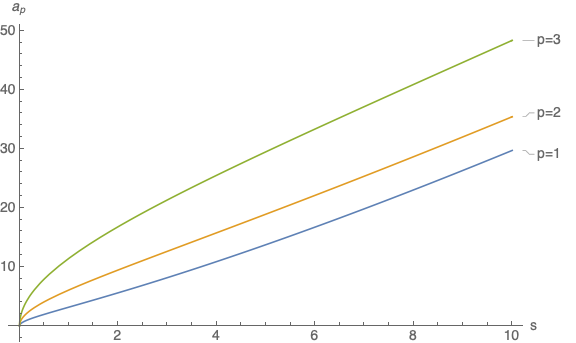}} \\
The graph of $a_p(s)$ for $p=1,2,3$. 
\end{center}

\begin{center}
\scalebox{0.6}{\includegraphics*{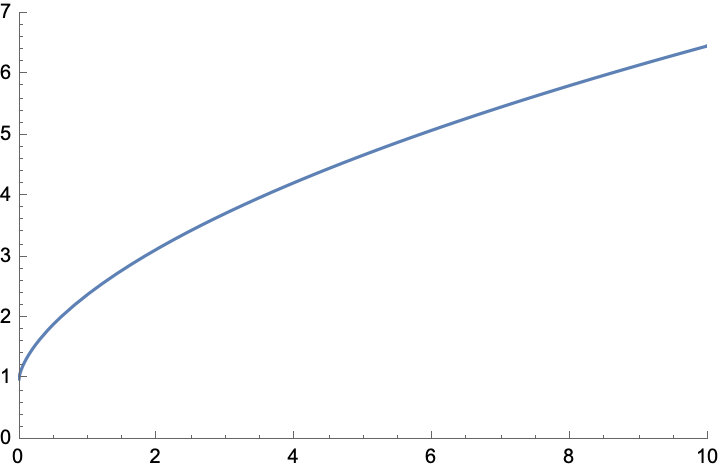}} \\
The graph of $m_p$ as a function of $p$. 
\end{center}

We next prove the  smoothness of $F$. 
\begin{lemma}\label{smoothness}
Suppose $F\in W^{1,2}_{loc}(\mP)$ satisfies (\ref{2.8})-(\ref{2.9}). Then $F$ is smooth.
\end{lemma}
\noindent{\bf Proof.} As we computed above, (\ref{2.8}) is an elliptic equation for $F$ (c.f. \cite[(15.1), pp.389]{AIM}), so it follows from the Caccioppoli-type estimates that $F\in W^{2,2}_{loc}(\ID )$. We next consider the function
\[
|F_z|^2=\eta^2a_p^2\left(\frac{F_\zbar-\sigma}{\eta}\right)=\eta^2\tilde{a}_p\left(\frac{F_\zbar-\sigma}{\eta}\right),
\]
where
\begin{eqnarray*}
\tilde{a}_p(s)&=&a^2_p(s)=(s+e^p)^2\left(\log^2(s+e^p)-p^2\right),\\
\tilde{a}_p'(s)&=&2(s+e^p)\left(\log^2(s+e^p)+\log(s+e^p)-p^2\right),\\
\tilde{a}_p''(s)&=&2\left(\log^2(s+e^p)+3\log(s+e^p)+1-p^2\right)\geq0.
\end{eqnarray*}
Then
\[
\min_{s\geq0}\tilde{a}_p^\prime(s)=\tilde{a}_p^\prime(0)=2pe^p>0,
\]
thus $\tilde{a}_p$ is invertible and we can write $\mathcal{B}_p=\tilde{a}_p^{-1}$, and then
\begin{equation}\label{2.10}
F_\zbar =\mathcal{B}_p\left(\frac{|F_z|^2}{\eta^2}\right)\eta+\sigma.
\end{equation}
Note here $\mathcal{B}_p(t^2)=\mathcal{A}_p(t)$, thus
\begin{equation}\label{2.11}
\mathcal{A}_p^\prime(t)=2t\mathcal{B}_p^\prime(t^2).
\end{equation}
As $F\in W^{2,2}_{loc}(\ID )$, we may differentiate both sides of (\ref{2.10}) by $x$, and get
\begin{equation}\label{2.12}
(F_x)_\zbar =\mathcal{B}_p^\prime\left(\frac{|F_z|^2}{\eta^2}\right)\frac{\overline{F_z}}{\eta}(F_x)_z+\mathcal{B}_p^\prime\left(\frac{|F_z|^2}{\eta^2}\right)\frac{F_z}{\eta}\overline{(F_x)_z}+\phi(z),
\end{equation}
where
\[
\phi(z)=\sigma_x+\mathcal{B}_p\left(\frac{|F_z|^2}{\eta^2}\right)\eta_x-2\mathcal{B}_p^\prime\left(\frac{|F_z|^2}{\eta^2}\right)\frac{|F_z|^2}{\eta^2}\eta_x\in W_{loc}^{1,2}(\mP).
\]
It follows from (\ref{2.11}) that
\[
\mathcal{B}_p^\prime\left(\frac{|F_z|^2}{\eta^2}\right)\frac{|F_z|}{\eta}+\mathcal{B}_p^\prime\left(\frac{|F_z|^2}{\eta^2}\right)\frac{|F_z|}{\eta}\leq\mathcal{A}_p^\prime\left(\frac{|F_z|}{\eta}\right)\leq k_p.
\]
Thus (\ref{2.12}) is again an elliptic equation for the function $F_x$, thus $F_x\in W^{2,2}_{loc}(\ID )$. The same argument applies on the function $F_y$, so $F\in W_{loc}^{3,2}(\ID )$. Thus we can differentiate (\ref{2.12}) again, and then we conclude that $F$ is smooth by a simple induction. \hfill $\Box$

\bigskip

By (\ref{2.7}) we can see that
\[
\mu_f=\frac{\overline{F_z}}{(F_\zbar-\sigma+e^p\eta)\left[\log\left(\frac{F_\zbar-\sigma}{\eta}+e^p\right)+p\right]},
\]
and so $\mu_f$ is also smooth. This implies $|\mu_f|<1$ and then $f$ is diffeomorphic, see \cite[Lemma 2.4 and Theorem 1]{MY1}.\\

We remark that the standing assumption $\IK(z,f)\exp(p\IK(z,f))\in L_{loc}^2(\mP)$ can be reduced. In fact, if $\exp(p\IK(z,f))\in L_{loc}^q(\mP)$ for some $q>1$, then we have $F\in W_{loc}^{1,q}(\mP)$. By (\ref{2.7}),
\begin{equation}\label{2.13}
\frac{F_\zbar-\sigma+e^p\eta}{|F_z|}=\frac{1-|\mu_f|^2}{2p|\mu_f|}.
\end{equation}
Let $k_0<1$ be sufficiently close to $1$ so that the right hand side of (\ref{2.13}) is smaller than $q-1$ if $|\mu_f|>k_0$; while for $|\mu_f|<k_0$, we have $\IK(z,f)\leq\frac{1+k_0^2}{1-k_0^2}$. Then we can write
\[
F_\zbar=\nu F_z+\tilde{\sigma},
\]
where $|\nu|<q-1$, $\tilde{\sigma}\in L^2(\mP)$. Thus the operator $I-\nu\mathcal{S}$ is invertible for $F_\zbar$, see \cite[Theorem 14.0.4]{AIM}. In fact, with the help of a cut-off function $\varphi\in C_0^\infty(\mP)$, $\varphi\equiv1$ in any given $\Omega\Subset\mP$, we have
\[
(\varphi F)_\zbar=(I-\nu\mathcal{S})^{-1}(\varphi_\zbar F-\nu\varphi_zF+\varphi\tilde{\sigma})\in L^2(\mP).
\]
This proves $F\in W_{loc}^{1,2}(\mP)$. This allows us to now conclude the following:
\begin{theorem}\label{diffeothm}
Let $f$ be a homeomorphic solution to (\ref{2.1}), and $\exp(p\IK(z,f))\in L_{loc}^q(\mP)$ for some $q>1$, then $f$ is diffeomorphic. In particular, any locally quasiconformal solution to (\ref{2.1}) is diffeomorphic.
\end{theorem}

\subsection{The $\mu$ equations.}
We now assume $f$ is a diffeomorphic minimiser of $\mE_p$. Then it satisfies equation (\ref{2.1}). This gives the following pointwise equation
\[
\left(\exp(p\IK(z,f))\frac{2p\overline{\mu_f}}{1-|\mu_f|^2}\eta\right)_\zbar=\left(\exp(p\IK(z,f))\right)_z\eta.
\]
Note this holds in every fundamental polyhedron $\mP$, thus  everywhere in $\ID$. This leads to an equation for $\mu=\mu_f$:
\begin{equation}\label{2.14}
\gamma(|\mu|)\mu_z=\alpha(|\mu|)\mubar\mu_\zbar-\beta(|\mu|)\mu^2\overline{\mu_\zbar}+\phi,
\end{equation}
where
\begin{eqnarray*}
\gamma(t)&=&1+(4p-3)t^2+(4p^2-4p+3)t^4-t^6,\\
\alpha(t)&=&(1-t^2)^3,\\
\beta(t)&=&2p(1+2pt^2-t^4),\\
\phi&=&-(1-|\mu|^2)^2(1+(2p-1)|\mu|^2)\mu\frac{\eta_z}{\eta}-(1-|\mu|^2)^3|\mu|^2\frac{\eta_\zbar}{\eta}.
\end{eqnarray*}
We consider
\[
A(t):=\frac{\alpha(t)t+\beta(t)t^2}{\gamma(t)}=\frac{t(1+(2p+1)t-t^2-t^3)}{1+t+(2p-1)t^2-t^3}.
\]
Then
\begin{equation}\label{2.15}
\frac{1+A^2}{1-A^2}\leq C\frac{1}{1-A}=C\frac{1+t+(2p-1)t^2-t^3}{(1-t^2)^2}.
\end{equation}
Thus we can rewrite (\ref{2.14}) as
\begin{equation}\label{2.16}
\mu_z=\nu\mu_\zbar+\phi,
\end{equation}
where $\phi\in L_{loc}^\infty(\ID)$, and
\begin{equation}\label{2.17}
\frac{1+|\nu|^2}{1-|\nu|^2}\leq C\IK(z,f)^2.
\end{equation}

On the ``$h$ side'', a diffeomorphic minimiser satisfies equation (\ref{2.2}), and then it satisfies Ahlfors-Hopf equation (\ref{2.3}), which holds in $\ID$. This subsequently implies that the Beltrami coefficient of $h$ has the following form:
\begin{equation}\label{2.18}
\mu_h=|\mu_h|\frac{\overline{\Phi}}{|\Phi|}.
\end{equation}
One can compare this and equation (\ref{teich1}) for the form of the Beltrami equation for  Teichm\"uller extremal mappings. Also, since $\Phi$ is holomorphic, (\ref{2.3}) implies
\[
\left(\exp(p\IK(w,h))h_w\overline{h_\wbar}\eta(h)\right)_\wbar=\Phi_\wbar=0.
\]
This leads to the equation
\begin{equation}\label{2.19}
h_{w\wbar}+(\log\lambda)_z(h)h_wh_\wbar=0,
\end{equation}
where
\begin{equation}\label{2.20}
\lambda(z)=\exp(p\IK(z,f))\eta(z).
\end{equation}
This is the {\it tension equation} for the harmonic mapping with respect to the metric $\lambda(z)|dz|^2$, \cite{M3}. We record these results in the following:
\begin{theorem}
A diffeomorphic minimiser $f$ of $\mE_p$ satisfies equation (\ref{2.14}), while its inverse $h$ satisfies equations (\ref{2.18}) and (\ref{2.19}). In particular, $h$ is a harmonic mapping from $\ID$ to $(\ID,\lambda)$, where $\lambda$ is the metric defined by (\ref{2.20}).
\end{theorem}

\section{Solution to the Ahlfors-Hopf equation}
In this section we discuss the holomorphic Ahlfors-Hopf differential
\begin{equation}\label{3.1}
\Phi=\exp(p\IK(w,h))h_w\overline{h_\wbar}\eta(h).
\end{equation}
We wish to prove the following:
\begin{theorem}\label{Hopfdiff}
Let $h:\oD\to\oD$ be a homeomorphism which satisfies equation (\ref{3.1}). Then $h$ is diffeomorphic in $\ID$.
\end{theorem}

This, together with Theorem \ref{AHdiff} and our previous work \cite{MY5} obtaining the uniqueness, completes the main theorem of this article:
\begin{theorem}\label{MainTheorem}
Let $R,S$ be analytically finite Riemann surfaces, $f_0:R\to S$ is a finite distortion homeomorphism, and
\[
\mE_p(f_0)=\int_R\exp(p\IK(z,f)) \; d\sigma_R(z)<\infty,\quad0<p<\infty.
\]
Then, in the homotopy class of $f_0$, there exists a unique homeomorphic minimiser $f$ of $\mE_p$. Moreover, the minimiser $f$ is a diffeomorphism.
\end{theorem}

We already know $h$ is a homeomorphism as it has degree $1$ and its inverse has exponentially integrable distortion.  A homeomorphism being a diffeomorphism is a local property,  we only need to prove that every point $w\in\ID$ has a neighbourhood $w\in\Omega\Subset\ID$ where $h$ is diffeomorphic.  On such a neighbourhood $\Phi$ is continuous and bounded. Because $\Phi$ is holomorphic, the zero set $Z:=\{w\in\ID:\Phi(w)=0\}$ is discrete in $\ID$.  Our argument will break down into two cases depending on whether $w\not\in Z$ or $w\in Z$.  In fact  $Z$ will be our `bad set' and this is because in general the extremal quasiconformal mapping is not a diffeomorphism on $Z$ -- the Beltrami coefficient is $k\overline{\psi}/|\psi|$, $0<k<1$ a real constant and so fails to be continuous where $\psi=0$. Thus one should not expect uniform gradient estimates on the distortion to be available as $p\to\infty$.  

We use the following strategy for our proof.

 \begin{enumerate}
 \item Choose $\Omega\Subset\ID$.  We have $h:\overline{\Omega} \to h(\overline{\Omega} )$ a homeomorphism satisfying (\ref{3.1}) on $\Omega$ and where $\Phi\in L^\infty(\Omega)$. Typically $\Omega$ or $h(\Omega)$ will be a disk and so both will be Jordan domains. The two cases are  $Z\cap \Omega=\emptyset$ and $Z\cap \Omega\neq\emptyset$.
\item In $\Omega$ we formulate an approximate extremal problem parameterised by $\lambda\in (0,1]$ whose variational equations are uniformly elliptic (with ellipticity bound depending on $\lambda$). We obtain our original problem and variational equations as $\lambda\to1$. 
\item For each $\lambda$ we use the nonlinear Riemann mapping theorem \cite[\S 9.2.1 \& 9.2.2]{AIM} with a two point normalisation determined by $h$  to find $\lambda$-solutions $h^\lambda$ to the elliptic variational equations.  These are guaranteed to be diffeomorphisms.
\item  We prove that as $\lambda\to1$ the functions $h^\lambda$ converge to some $H$ which is also diffeomorphic in $\Omega$, and which has the same Ahlfors-Hopf differential (\ref{3.1}) as $h$.  This establishes $H$ as a diffeomorphic solution to  (\ref{3.1}) in $\Omega$.
\item We finally prove a uniqueness theorem to show that $h|\Omega=H$, so $h$ is diffeomorphic in $\Omega$.  
\end{enumerate}
We will use different methods for points in $\ID\setminus Z$ and $Z$.

\subsection{The $\lambda$-equations.}
We will continue to use $f$ and $h$ for mappings in these variational problems,  though the  hypotheses of Theorem \ref{Hopfdiff} provide us with $h$ and $f=h^{-1}$ defined on $\ID$ we set them aside and just reference them as `boundary data' or in the mapping problem we ultimately consider,  as defining the domains in question. Note that these mappings defined on $\ID$ provide a barrier for the local problem (with the boundary values they define) and so the minimum value of the functional we describe is certainly finite.

From now on, $\Omega,\Omega'\Subset\ID$ will always be planar Jordan domains. We consider the following $\lambda$-problems and their inverses. We consider the problem of minimising the following functional among all Sobolev homeomorphisms $f:\Omega\to \Omega'$ with $h=f^{-1}:\Omega'\to\Omega$ with fixed given boundary data $f|_{\partial\Omega}$.   

\begin{equation}\label{3.2}
\int_\Omega\exp\left(p\frac{\lambda^2+|\mu_f|^2}{\lambda^2-|\mu_f|^2}\right)\eta(z) \; dz, \;{\rm and} \;      \int_{\Omega'}\exp\left(p\frac{\lambda^2+|\mu_h|^2}{\lambda^2-|\mu_h|^2}\right)J(w,h)\eta(h) \; dw
\end{equation}
 Inner variation of each of these functionals gives the following Euler-Lagrange equations.
\begin{eqnarray}\label{3.3}
\lefteqn{\int_\Omega\exp\left(p\frac{\lambda^2+|\mu_f|^2}{\lambda^2-|\mu_f|^2}\right)(\eta\varphi)_z \; dz}\nonumber \\ &=&\int_\Omega\exp\left(p\frac{\lambda^2+|\mu_f|^2}{\lambda^2-|\mu_f|^2}\right)\frac{2p\lambda^2\overline{\mu_f}(1-|\mu_f|^2)}{(\lambda^2-|\mu_f|^2)^2}\eta\varphi_\zbar \; dz,\label{3.3}
\end{eqnarray}
\begin{equation}\label{3.4}
\int_{\Omega'}\exp\left(p\frac{\lambda^2+|\mu_h|^2}{\lambda^2-|\mu_h|^2}\right)\left(\frac{1-|\mu_h|^2}{\lambda^2-|\mu_h|^2}\right)^2h_w\overline{h_\wbar}\eta(h)\psi_\wbar \; dw=0,
\end{equation}
for any $\varphi\in C_0^\infty(\Omega)$ and $\psi\in C_0^\infty(\Omega')$, respectively. Note a minimiser of (\ref{3.2}) is not automatically variational. However, for any given holomorphic mapping $\Phi$ defined in $\Omega'$, we may consider the equation
\begin{equation}\label{3.5}
\Phi=\exp\left(p\frac{\lambda^2+|\mu_h|^2}{\lambda^2-|\mu_h|^2}\right)\left(\frac{1-|\mu_h|^2}{\lambda^2-|\mu_h|^2}\right)^2h_w\overline{h_\wbar}\eta(h).
\end{equation}
We notice that (\ref{3.5}) defines an implicit functional relationship between $h_w$ and $h_\wbar$. We now turn to study this function.

\subsection{The non-linear Beltrami operators $\mathcal{B}^\lambda$, $0<\lambda \leq 1$.}
To analyse the equation (\ref{3.5}) we first take the absolute values and consider
\begin{equation}\label{3.6}
\exp\left(p\frac{\lambda^2x^2+y^2}{\lambda^2x^2-y^2}\right)\left(\frac{x^2-y^2}{\lambda^2x^2-y^2}\right)^2xy=k,
\end{equation}for $0\leq y\leq \lambda x$ as a function of $x\geq 0$ and $k\geq0$. Notice that for fixed $x>0$ and $k>0$, the left hand side of (\ref{3.6}) is strictly increasing with $y$ from $0$ to $+\infty$, therefore there is a unique $y=A_k^\lambda(x)$ solving (\ref{3.6}).  The implicit function theorem assures us that the assignment $y=\mathcal{A}_k^\lambda(x)$ is a smooth function in all variables on the positive quadrant
\[
Q^+=\{(x,k)\in  (0,\infty)\times (0,\infty)\}.
\]
We define
\begin{equation}\label{3.7}
v=v_k^\lambda(x) = \frac{\mathcal{A}_k^\lambda(x)}{x}, \hskip20pt (x,k)\in  (0,\infty)\times (0,\infty).
\end{equation}
Then $v$,  $0<v<\lambda$,  is a $C^\infty$ function in $x$ and $k$ for all $\lambda$ on $Q^+$. Further on $Q^+$ the function $v$ also solves the functional equation
\begin{equation}\label{3.8}
\exp\left(p\frac{\lambda^2+v^2}{\lambda^2-v^2}\right)\left(\frac{1-v^2}{\lambda^2-v^2}\right)^2v= \frac{k}{x^2}.
\end{equation}
It is immediate that for $x>0$ the assignment $v_0^\lambda(x)=0$ is the continuous extension of $v$ to $(x,k)\in (0,\infty)\times [0,\infty)$.  Further,  if $k\neq 0$ it is clear that the assignment $v_k^\lambda(0)=\lambda$ is also continuous.  So defined we have that $v_k^\lambda(x) $ is continuous on $\overline{Q^+}\setminus \{(0,0)\}$.\\

We implictly differentiate the equation (\ref{3.8}) to obtain the following two equations.
\begin{eqnarray}\label{3.9}
\frac{\partial v}{\partial x} \; \Big[\frac{4p\lambda^2v}{(\lambda^2-v^2)^2}+\frac{4v(1-\lambda^2)}{(\lambda^2-v^2)(1-v^2)}+\frac{1}{v}\Big]=-\frac{2}{x}, \\
\frac{\partial v}{\partial k} \; \Big[\frac{4p\lambda^2v}{(\lambda^2-v^2)^2}+\frac{4v(1-\lambda^2)}{(\lambda^2-v^2)(1-v^2)}+\frac{1}{v}\Big]= \frac{1}{k}. \label{3.10}
\end{eqnarray}
These are both done in the same way. For each $k>0$, note that by (\ref{3.8}), $v=\lambda$ when $x=0$, and $v=0$ when $x=\infty$. Combine (\ref{3.8}) and (\ref{3.9}) we have
\[
\frac{\partial v}{\partial x} \; \Big[\frac{4p\lambda^2v}{(\lambda^2-v^2)^2}+\frac{4v(1-\lambda^2)}{(\lambda^2-v^2)(1-v^2)}+\frac{1}{v}\Big]=-2\sqrt{\frac{v}{k}\exp\left(p\frac{\lambda^2+v^2}{\lambda^2-v^2}\right)}\left(\frac{1-v^2}{\lambda^2-v^2}\right)
\]
So we have $\frac{\partial v}{\partial x}=-\infty$ when $x=0$ and $\frac{\partial v}{\partial x}=0$ when $x=\infty$. Also note
\[
\frac{4p\lambda^2v}{(\lambda^2-v^2)^2}+\frac{4v(1-\lambda^2)}{(\lambda^2-v^2)(1-v^2)}+\frac{1}{v}\geq\lambda^2\Big[\frac{4pv}{(\lambda^2-v^2)^2}+\frac{1}{v}\Big]\geq\min\left\{p\lambda^2,4\lambda^2\right\}.
\]
Thus by (\ref{3.10}) we have
\begin{equation}\label{3.11}
\frac{\partial v}{\partial k}\leq\frac{c_p}{k},\quad c_p=\max\left\{\frac{1}{p\lambda^2},\frac{1}{4\lambda^2}\right\}.
\end{equation}
Combine (\ref{3.8}) and (\ref{3.10}) we also get
\[
\frac{\partial v}{\partial k} \; \Big[\frac{4p\lambda^2v}{(\lambda^2-v^2)^2}+\frac{4v(1-\lambda^2)}{(\lambda^2-v^2)(1-v^2)}+\frac{1}{v}\Big]\exp\left(p\frac{\lambda^2+v^2}{\lambda^2-v^2}\right)\left(\frac{1-v^2}{\lambda^2-v^2}\right)^2v= \frac{1}{x^2}.
\]
For fixed $x>0$, when $k=0$, $v=0$, so $\frac{\partial v}{\partial k}=\frac{\lambda^2}{e^px^2}$; when $k=\infty$, $v=\lambda$, so $\frac{\partial v}{\partial k}=0$.\\

We now put together the facts we have established about $v$.

\begin{lemma}\label{vlemma}
Let $v=v^\lambda_k(x)$. Then 
\begin{enumerate}
\item  $v\in C^{\infty}(Q^+)$,  $0\leq v \leq \lambda$ and on $Q^+$,  $v$ is strictly decreasing in $x$ and strictly increasing in $k$.
\item For fixed $k>0$, $v\to \lambda$, $\frac{\partial v }{\partial x} \to -\infty$ as $x\to0$; $v\to0$, $\frac{\partial v }{\partial x} \to 0$ as $x\to\infty$.
\item For fixed $x>0$,  $v\to 0$, $\frac{\partial v }{\partial k}\to\frac{\lambda^2}{e^px^2}$ as $k\to 0$; $v(x) \to \lambda$, $\frac{\partial v }{\partial k}\to0$ as $k\to \infty$.
\end{enumerate}
\end{lemma}

\medskip

Next,  given a holomorphic function $\Phi$ in $\Omega'$ we define the functions $\mathcal{B}^\lambda:\Omega'\times\IC\times\IC\to \IC$ by
\begin{equation}\label{3.12}
\mathcal{B}^\lambda(w,\tau,\xi) = \frac{\overline{\Phi(w)}}{|\Phi(w)|}  v^\lambda_{\frac{|\Phi(w)|}{\eta(\tau)}}(|\xi|) \xi.
\end{equation}

We need to establish certain properties of the function $\mathcal{B}^\lambda$.  Of course it is immediate that $0\leq \mathcal{B}^\lambda \leq \lambda$ and that $\mathcal{B}^\lambda$ is continuous on $\Omega'\times \IC\to \IC$.

\subsection{$\mathcal{B}^\lambda$ is Lipschitz in $\xi\in\IC$.}
We first compute 
\begin{eqnarray}\label{3.13}
\frac{|\mathcal{B}^\lambda(w,\tau,\zeta)-\mathcal{B}^\lambda(w,\tau,\xi)| }{|\zeta-\xi|}& = & \frac{ \big|v(|\zeta|)  \zeta -v(|\xi|)  \xi \big|}{|\zeta-\xi|}.
\end{eqnarray}
 We put  $|\zeta|=t$,  $|\xi|=s$, $\alpha=v(t)$ and $\beta=v(s)$.  Then there is a $\theta\in[0,2\pi]$ such that
$\zeta\cdot\xi= st\cos(\theta),$
and
\begin{eqnarray*}
\frac{ \big|v(|\zeta|)  \zeta -v(|\xi|)  \xi \big|^2}{|\zeta-\xi|^2} & = & \frac{\alpha^2 t^2+\beta ^2 s^2 - 2 \alpha\beta st \cos(\theta)}{t^2+s^2-2st\cos(\theta)}:=F(\theta).
\end{eqnarray*}
We differentiate this with respect to $\theta$ to see
\[
\frac{d}{d\theta}F(\theta)=\frac{2  st(\alpha -\beta )(s^2\beta -t^2\alpha )}{(t^2+s^2-2st\cos(\theta))^2}\sin(\theta).
\]
At this point we claim that
\begin{equation}\label{3.14}
(\alpha -\beta )(s^2\beta -t^2\alpha )\geq0.
\end{equation}
To see this define $w(x)=x^2v(x) \leq \lambda x^2$. Then differentiating the functional equation shows that $w$ is increasing,  $w'(x)=$
\[ {  \frac{8x^3w^3[p\lambda^2(x^4-w^2)+(1-\lambda^2)(\lambda^2x^4-w^2)]}{4p\lambda^2x^4w^2(x^4-w^2)+4w^2(1-\lambda^2)x^4(\lambda^2x^4-w^2)+(\lambda^2x^4-w^2)^2(x^4-w^2)}>0}. \]
Now assume $t\leq s$, then $\alpha = v(t)\geq v(s)=\beta$, 
$s^2\beta=w(s)\geq w(t)=t^2\alpha$,
and vice versa. So (\ref{3.14}) follows. 
Assume $\zeta\neq\xi$.  Then
\[
\frac{d}{d\theta}\frac{\big|v(|\zeta|)\zeta-v(|\xi|)\xi\big|^2}{|\zeta-\xi|^2}=G(|\zeta|,|\xi|,\cos(\theta))\sin\theta,
\]
where $G(|\zeta|,|\xi|,\cos(\theta))$ is always non-negative. In the period $\theta\in[0,2\pi]$ we see that $F(\theta)$ is increasing when $\theta\in [0,\pi]$ and decreasing on $[\pi,2\pi]$ so we get the maximum of $F(\theta)$ at $\theta=\pi$. In particular we can now write (\ref{3.13}) as
\begin{equation}\label{3.15}
\frac{|\mathcal{B}^\lambda(w,\tau,\zeta)-\mathcal{B}^\lambda(w,\tau,\xi)| }{|\zeta-\xi|} \leq\frac{ v(t) t +v(s)s }{t+s} \leq\max\{v(|\zeta|),v(|\xi|)\},
\end{equation}
whenever $\zeta\neq\xi$.  There is actually more useful information here that we will use to establish uniqueness later.  Consider the function $u(s,t)=\frac{ v(t) t +v(s)s }{t+s}$,  smooth on $s,t>0$.  We compute
\begin{eqnarray*}
u_s(s,t) & = & \frac{s (s+t) v'(s)+t (v(s)-v(t))}{(s+t)^2}.
\end{eqnarray*}
Since $v'(s)<0$ we have $u_s(s,t) =0$ only if $v(s)\geq v(t)$,  thus $s\leq t$.  By symmetry $\nabla u = 0$ implies that $s=t$.  We wish to find the maximum value of $u$ on $[t_0,\infty)\times (0,\infty)$,  where $u$ is also continuous.  If this value occurs at an interior critical point,  then we have $s=t\geq t_0$ and that value is $v(t_0)$. 
Otherwise the maximum occurs on the boundary $s=0$ or $t=t_0$.  When $s=0$,  $u(t,s)=v(t)\leq v(t_0)$. We now need to consider $u_s(s,t_0)$, differentiable for $s>0$.  At the boundary $s=0$ or $s=+\infty$ we again recover the value $v(t_0)$ (since $v\to 0$ as $s\to\infty$).  It now follows that the minimum value occurs at a point in the interval $(0,t_0]$.  Call this point $s_0$.  If $v(s_0)\geq \frac{1}{2}(\lambda + v(t_0))$,  then 
\begin{eqnarray*} \frac{ v(t_0) t_0 +v(s_0)s_0 }{t_0+s_0} &=&v(s_0) +  \frac{t_0}{t_0+s_0}(v(t_0) - v(s_0)) \leq  v(s_0) +  \frac{t_0}{t_0+s_0}\;\frac{v(t_0)-\lambda}{2} \\
&\leq&\frac{1}{2}(\lambda + v(t_0)).
\end{eqnarray*}
Otherwise 
\[
\frac{ v(t_0) t_0 +v(s_0)s_0 }{t_0+s_0} \leq \max\{v(t_0),v(s_0)\}=v(s_0)\leq \frac{1}{2}(\lambda + v(t_0)).
\]
In either case we deduce the following lemma.
\begin{lemma}\label{qrlemma} For $w\in\Omega$ and $\zeta,\xi\in \IC$  and $\lambda\in (0,1]$ we have
\begin{equation}\label{3.16}
\frac{|\mathcal{B}^\lambda(w,\tau,\zeta)-\mathcal{B}^\lambda(w,\tau,\xi)| }{|\zeta-\xi|}   \leq \max\{v(|\zeta|),v(|\xi|)\} \leq \lambda.
\end{equation}
Moreover,  if  $|\xi|\geq t_0>0$,  then for $w\in \Omega$
\begin{equation}\label{3.17}
\frac{|\mathcal{B}^\lambda(w,\tau,\zeta)-\mathcal{B}^\lambda(w,\tau,\xi)| }{|\zeta-\xi|}   \leq \frac{1}{2}(1+ v(t_0) )  < 1.
\end{equation}
The last estimate is locally uniform in $w$,  and  uniform if $\Phi$ is bounded.
\end{lemma}
\noindent {\bf Proof.} In light of what we already have computed we only need to make the observation that $v^\lambda_k(x)$ is increasing in $k=|\Phi(w)|$ and that differentiating the functional equation shows that $v^\lambda_k(x)$ is strictly increasing in $\lambda$. \hfill $\Box$

\bigskip

Lemma \ref{qrlemma} shows that equation (\ref{3.5}) is an elliptic equation. Now the following follows from \cite[Theorem 9.0.3]{AIM}:
\begin{theorem}\label{lambdaExistence}
Let $0<\lambda<1$, $w_0\in\Omega'$, $z_0\in\Omega$,  $a\in\partial\Omega'$, $b\in\partial\Omega$. Then there exists a quasiconformal solution $h^\lambda:\overline{\Omega'}\to\overline{\Omega}$ to equation (\ref{3.5}) such that $h^\lambda(w_0)=z_0$, $h^\lambda(a)=b$.
\end{theorem}
The quasiconformality bound here depends on $\lambda$ and so we will look for uniform bounds elsewhere.

\subsection{Convergence of $h^\lambda$.}
For each $\lambda<1$, we let $h^\lambda$ be a solution to (\ref{3.5}) as in Theorem \ref{lambdaExistence}. In this section we prove:
\begin{lemma}\label{UniConv}
As $\lambda\to1$, $h^\lambda\to H$ uniformly in $\overline{\Omega'}$, where $H:\overline{\Omega'}\to\overline{\Omega}$ is also a homeomorphism, and $H(z_0)=w_0$, $H(a)=b$. Furthermore, their inverses $f^\lambda=(h^\lambda)^{-1}$ also converge uniformly in $\overline{\Omega}$ to $g=H^{-1}$.
\end{lemma}
Write
\[
\IK^\lambda(z,\cdot)=\frac{\lambda^2+|\mu(z,\cdot)|^2}{\lambda^2-|\mu(z,\cdot)|^2}.
\]
We claim that $h^\lambda$ have a uniform $W^{1,2}(\Omega')$ norm. In fact,
\begin{eqnarray*}
&&\int_{\Omega'}\|Dh^\lambda(w)\|^2 \; dw\;=\;\int_{\Omega'}\IK(w,h^\lambda) \; J(w,h^\lambda)dw\\
&\leq&\int_{\Omega'}\exp\left(p\IK^\lambda(w,h^\lambda)\right) \; J(w,h^\lambda)dw+C\\
&\leq&\int_{\Omega'\cap\{|\mu_{h^\lambda}|\geq\frac{\lambda}{2}\}}\frac{|\Phi(w)|}{|\eta(h^\lambda)|}\frac{\lambda^2-|\mu_{h^\lambda}|^2}{|\mu_{h^\lambda}|}\;dw \\ && +\int_{\Omega'\cap\{|\mu_{h^\lambda}|\leq\frac{\lambda}{2}\}}\exp\left(p\IK^\lambda(w,h^\lambda)\right) \; J(w,h^\lambda)dw+C\\
&\leq&\frac{3}{2}\pi\left\|\frac{1}{\eta}\right\|_\infty+\exp\left(\frac{5p}{3}\right)\pi+C<\infty.
\end{eqnarray*}
Let $f^\lambda=(h^\lambda)^{-1}$. Then
\begin{eqnarray*}
\int_\Omega\exp\left(p\IK(z,f^\lambda)\right) \; dz&\leq&\int_\Omega\exp\left(p\IK^\lambda(z,f^\lambda)\right) \; dz\\
&=&\int_\Omega\exp\left(p\IK^\lambda(w,h^\lambda)\right) \; J(w,h^\lambda)dw,
\end{eqnarray*}
which are also uniformly bounded as above. 

Thus $f^\lambda$ are bounded in $W^{1,Q}(\Omega)$ for $Q(t)=\frac{t^2}{\log(e+t)}$, and then (up to a subsequence) both $f^\lambda$ and $h^\lambda$ converge locally uniformly (this follows from Theorem \ref{thm4} and the related discussion in \S \ref{expsec}). We call the limit functions $g$ and $H$, respectively. By uniform convergence, they are continuous in $\Omega$ and $\Omega'$, respectively.\\

Next we wish to extend the functions $f^\lambda$ and $h^\lambda$ to some slightly larger domains so as to get uniform convergence in $\overline{\Omega}$ and $\overline{\Omega'}$. Typically this is by reflection. However, we cannot reflect over $\Omega$ and $\Omega'$, so we transform them to disks keeping the ellipticity estimates. 

Consider the conformal mappings $\phi:\oD\to\overline{\Omega}$ and $\psi:\overline{\Omega'}\to\oD$. Then each $\psi\circ f^\lambda\circ \phi$ and $\phi^{-1}\circ h^\lambda\circ\psi^{-1}$ are self-homeomorphisms of $\oD$, and for any $q<p/3$,
\begin{eqnarray*}
\lefteqn{\int_\ID\exp\left(q\IK(z,\psi\circ f^\lambda\circ\phi)\right) \; dz}\\
& = &\int_\ID\exp\left(q\IK(\phi,f^\lambda)\right) \; dz\\
&=&\int_\Omega\exp\left(q\IK(z,f^\lambda)\right)|D\phi^{-1}(z)|^2 \; dz\\
&\leq&\left(\int_\Omega\exp\left(3q\IK(z,f^\lambda)\right)\; dz\right)^{1/3}\left(\int_\Omega|D\phi^{-1}(z)|^3 \; dz\right)^{2/3}<\infty,
\end{eqnarray*}
and thus
\begin{eqnarray*}
\int_\ID|D\phi^{-1}\circ h^\lambda\circ\psi^{-1}(w)|^2 \; dw&=&\int_\Omega\IK(z,\psi\circ f^\lambda\circ\phi) \; dz\\
&\leq&\int_\Omega\exp\left(q\IK(z,\psi\circ f^\lambda\circ\phi)\right) \; dz+C<\infty.
\end{eqnarray*}
We remark that the integrability of $|D\phi^{-1}|^3$ arises from the famous conjecture of Brennan \cite{Brennan}. He proved that in order for $\int_\Omega|D\phi^{-1}|^{p_0}<\infty$, the largest possible value $p_0>3$ and then conjectured that it should always hold true for any $p_0<4$. 

Now by the same argument as in the proof of Corollary \ref{Extension}, we have $\psi\circ f^\lambda\circ \phi$ and $\phi^{-1}\circ h^\lambda\circ\psi^{-1}$ both converge uniformly in $\oD$, to some $\tilde{g}$ and $\tilde{H}$, respectively. 

Hence $g=\psi^{-1}\circ\tilde{g}\circ\phi^{-1}$ and $H=\phi\circ\tilde{H}\circ\psi$ satisfy the requirements of Lemma \ref{UniConv}.

\subsection{Convergence of derivatives.}

We continue with the notation established above. Notice that the following lemma requires $\Phi\in L^1(\Omega')$ according to \cite{MY3}, although we in fact have $\Phi$ bounded  in $\Omega'$.
\begin{lemma}
$h^\lambda$ is the unique minimiser for the inverse boundary value problem (\ref{3.2}) with its own boundary values.
\end{lemma}
\noindent{\bf Proof.}  This is a slight variation of our previous work \cite{MY3} so we do not give all details. It is worth noting that 
\[
\IK^\lambda(w,h)=\frac{\lambda^2+|\mu_h|^2}{\lambda^2-|\mu_h|^2}=\frac{(1-\lambda^2)-(1+\lambda^2)\IK(w,h)}{-(1+\lambda^2)+(1-\lambda^2)\IK(w,h)}
\]
is a convex function of $\IK(w,h)$, if $1\leq\IK(w,h)\leq\frac{1+\lambda^2}{1-\lambda^2}$. Thus the $\lambda$ problem (\ref{3.2}) can be regarded as a particular case of \cite[Theorem 8]{MY3} with $\Psi(\IK(w,h))=\exp(p\IK^\lambda(w,h))$.\hfill$\Box$

\bigskip

Now since $h^\lambda$ is a minimiser for the inverse $\lambda$ problem, its inverse $f^\lambda$ is also a minimiser for the $\lambda$ problem. However, to prove that $f^\lambda$ is inner variational, we need it to be in a slightly better space.

\begin{lemma}
There is a $q>1$ depending on $\lambda$ such that $\exp(p\IK^\lambda(z,f^\lambda))\in L_{loc}^q(\Omega)$.
\end{lemma}
\noindent{\bf Proof.} We can transform (\ref{3.5}) to the ${f^\lambda}$ side to obtain
\[
\exp(p\IK^\lambda(z,f^\lambda))\frac{|\mu_{f^\lambda}|(1-|\mu_{f^\lambda}|^2)}{(\lambda^2-|\mu_{f^\lambda}|)^2}\frac{\eta}{J(z,f^\lambda)}=|\Phi({f^\lambda})|\in L_{loc}^\infty(\Omega).
\]
Thus for $q>1$, and any $\tilde{\Omega}\Subset\Omega$,
\begin{eqnarray*}
\lefteqn{\int_{\tilde{\Omega}}\exp\left(pq\IK^\lambda(z,f^\lambda)\right) \; dz}\\
&\leq&\int_{{\tilde{\Omega}}\cap\{|\mu_{f^\lambda}|\leq\lambda/2\}}\exp\left(pq\IK^\lambda(z,f^\lambda)\right) \; dz+\int_{{\tilde{\Omega}}\cap\{|\mu_{f^\lambda}|>\lambda/2\}}\exp\left(pq\IK^\lambda(z,f^\lambda)\right) \; dz\\
&\leq& C_1+C_2\int_{\tilde{\Omega}}|\Phi({f^\lambda})|^q(J(z, f^\lambda))^q \; dz\leq C_1+C_2\int_{\tilde{\Omega}} (J(z,f^\lambda))^q \; dz.
\end{eqnarray*}
Since $f^\lambda$ is $\frac{1+\lambda}{1-\lambda}$-quasiconformal, $J(z,f^\lambda)\in L_{loc}^q(\Omega)$ for some $q>1$ depending on $\lambda$, see \cite[Corollary 13.2.5]{AIM}.\hfill$\Box$

\bigskip

The next lemma follows from exactly the same reasoning as for the exponential problem in \S 2.3. 
\begin{lemma}
Let $f^\lambda$ be a minimiser for the boundary value problem (\ref{3.2}), and $\exp(p\IK^\lambda(z,f^\lambda))\in L_{loc}^q(\Omega)$ for some $q>1$, then $f^\lambda$ satisfies the inner variational equation (\ref{3.3}).
\end{lemma}

\bigskip

We now have that $f^\lambda$ satisfies (\ref{3.3}) for every $\varphi\in C_0^\infty(\Omega)$. We rewrite this as
\begin{eqnarray*}
\lefteqn{\int_\Omega\exp\left(p\IK^\lambda(z,f^\lambda)\right)\frac{2p\lambda^2\overline{\mu_f}(1-|\mu_f|^2)}{(\lambda^2-|\mu_f|^2)^2}\eta\varphi_\zbar \; dz}\\ &=& \int_\Omega\left(\left(\exp\left(p\IK^\lambda(z,f^\lambda)\right)-e^p\right)\eta+\sigma^\lambda\right)\varphi_z \; dz,
\end{eqnarray*}
where
\[
\sigma^\lambda=-\mathcal{C}^*\left(\left(\exp\left(p\IK^\lambda(z,f^\lambda)\right)-e^p\right)\eta_z\right)\in L^2(\Omega).
\]
This gives an $F^\lambda\in W_{loc}^{1,q}(\Omega)$ such that
\begin{equation}\label{3.18}
F^\lambda_z=\exp\left(p\IK^\lambda(z,f^\lambda)\right)\frac{2p\lambda^2\overline{\mu_{f^\lambda}}(1-|\mu_{f^\lambda}|^2)}{(\lambda^2-|\mu_{f^\lambda}|^2)^2}\eta,
\end{equation}
\begin{equation}\label{3.19}
F^\lambda_\zbar=\left(\exp\left(p\IK^\lambda(z,f^\lambda)\right)-e^p\right)\eta+\sigma^\lambda.
\end{equation}
By (\ref{3.18}) and (\ref{3.19}),
\[
\left|\frac{F^\lambda_\zbar-\sigma^\lambda+e^p\eta}{F^\lambda_z}\right|=\frac{(\lambda^2-|\mu_{f^\lambda}|^2)^2}{2p\lambda^2|\mu_{f^\lambda}|(1-|\mu_{f^\lambda}|^2)}.
\]
If we choose $\lambda_0$ sufficiently close to $\lambda$, $\frac{(\lambda^2-|\mu_{f^\lambda}|^2)^2}{2p\lambda^2|\mu_{f^\lambda}|(1-|\mu_{f^\lambda}|^2)}<q-1$ if $\lambda_0<|\mu_{f^\lambda}|\leq\lambda$, and $\exp\left(p\IK^\lambda(z,f^\lambda)\right)\leq\exp\left(p\frac{\lambda^2+\lambda_0^2}{\lambda^2-\lambda_0^2}\right)$ if $0\leq|\mu_{f^\lambda}|\leq\lambda_0$. So we can write
\[
F^\lambda_\zbar=\kappa^\lambda F^\lambda_z+M^\lambda,
\]
where $F^\lambda\in W_{loc}^{1,q}(\Omega)$, $|\kappa^\lambda|<q-1$, $M^\lambda\in L^2(\Omega)$. This implies $F^\lambda\in W_{loc}^{1,2}(\Omega)$. Again by the same reasoning as we used in \S 2.3, each $F^\lambda$ is smooth.\\

Our next aim is to give a uniform estimate for $F^\lambda$ as $\lambda\to1$. From (\ref{3.18}) and (\ref{3.19}),
\begin{equation}\label{3.20}
F^\lambda_\zbar=\mathcal{A}^\lambda_p\left(\left|\frac{F^\lambda_z}{\eta}\right|\right)\eta+\sigma^\lambda,
\end{equation}
where $\mathcal{A}^\lambda_p$ is the inverse of the increasing function of $s\geq0$,
\begin{equation}\label{3.21}
a_p^\lambda(s)=(s+e^p)\sqrt{\log^2(s+e^p)-p^2}\frac{(1-\lambda^2)\log(s+e^p)+(1+\lambda^2)p}{2p\lambda},
\end{equation} 

$\mathcal{A}^\lambda$ satisfies $\mathcal{A}^\lambda(0)=0$ and $0\leq(\mathcal{A}_p^\lambda)'\leq k_p<1$, where $k_p$ is uniform as $\lambda\to1$. Thus we can write $\nu^\lambda=\mathcal{A}_p^\lambda\left(\left|\frac{F_z^\lambda}{\eta}\right|\right)\frac{\eta}{F_z^\lambda}$, and then
\[
F^\lambda_\zbar = \nu^\lambda F^\lambda_z+\sigma^\lambda,
\]
where
\[
|\nu^\lambda|=\left|\frac{\mathcal{A}_p^\lambda\left(\left|\frac{F_z}{\eta}\right|\right)-\mathcal{A}_p^\lambda(0)}{\frac{F_z}{\eta}-0}\right|\leq\sup(\mathcal{A}_p^\lambda)'\leq k_p<1,
\]
and $\sigma^\lambda\in L^2(\Omega)$ is also uniform for $\lambda$. Let $\phi\in C_0^\infty(\Omega)$ be a cut-off function, then
\[
(\phi F^\lambda)_\zbar = \nu^\lambda(\phi F^\lambda)_z+\phi_\zbar F^\lambda-\nu^\lambda\phi_z F^\lambda+\phi\sigma^\lambda.
\]
Again we can write
\[
\psi^\lambda=\phi_\zbar F^\lambda-\nu^\lambda\phi_z F^\lambda+\phi\sigma^\lambda\in L^2(\IC).
\]
Then
\[
(\phi F^\lambda)_\zbar=(\mathbf{I}-\nu^\lambda\mathcal{S})^{-1}\psi^\lambda\in L^2(\IC).
\]
This proves that the $L^2_{loc}(\Omega)$ norms of $F^\lambda_z,F^\lambda_\zbar$ are independent of $\lambda$. We next differentiate both sides of (\ref{3.20}). To do this we define $\mathcal{B}_p^\lambda$ such that $\mathcal{B}_p^\lambda(t^2)=\mathcal{A}_p^\lambda(t)$. Then
\begin{eqnarray*}
(F_x^\lambda)_\zbar &=&(\mathcal{B}_p^\lambda)'\left(\frac{|F_z|^2}{\eta^2}\right)\frac{\overline{F^\lambda_z}}{\eta}(F^\lambda_x)_z+(\mathcal{B}_p^\lambda)'\left(\frac{|F^\lambda_z|^2}{\eta}\right)\frac{F^\lambda_z}{\eta^2}\overline{(F^\lambda_x)_z}\\
&&+\mathcal{B}_p^\lambda\left(\frac{|F_z|^2}{\eta^2}\right)\eta_x-2(\mathcal{B}_p^\lambda)'\left(\frac{|F^\lambda_z|^2}{\eta^2}\right)\frac{|F^\lambda_z|^2}{\eta^2}\eta_x+\sigma^\lambda_x.
\end{eqnarray*}
Again we write this in the form
\[
(F^\lambda_x)_\zbar=\tilde{\nu}^\lambda(F^\lambda_x)_z+\tilde{\sigma}^\lambda,
\]
which is again an elliptic equation and following the same argument as above we obtain a uniform estimate on $(F^\lambda_x)_z,(F^\lambda_x)_\zbar\in L^2_{loc}(\Omega)$. This same argument works for $F^\lambda_y$, and then inductively and we see that we can differentiate the equation infinitely many times to get uniform $W_{loc}^{k,p}(\Omega)$ norms of $F^\lambda$, for all $k\geq1$ and $p\geq1$. This gives the following:
\begin{lemma}
For all $\lambda<1$, $F^\lambda$ is smooth and has a locally uniform bound for $C^k(\Omega)$, for all $k=1,2,\dots$.
\end{lemma}

\bigskip

We now have, up to a subsequence, $F^\lambda\to F$ locally uniformly, and $F$ is smooth. In particular, $\mu_{f^\lambda}$ also converge to some smooth $\mu$. We also have $f^\lambda\to g$ in $W^{1,P}(\ID)$. From the ``good approximation lemma'' (see \cite[Lemma 5.3.1]{AIM}), $\mu_g=\mu$. This proves that $g$ is a diffeomorphism, and then so is $H=g^{-1}$.\\

We have proved the following:
\begin{theorem}\label{diffeoExistence}
Assume $\Phi\in L^1(\Omega')$. Let $w_0\in\Omega'$, $z_0\in\Omega$,  $a\in\partial\Omega'$, $b\in\partial\Omega$. Then there exists a homeomorphic solution $H:\overline{\Omega'}\to\overline{\Omega}$ to the equation defining a holomorphic Ahlfors-Hopf differential (\ref{3.1}), such that $H(w_0)=z_0$, $H(a)=b$, and $H$ is diffeomorphic in $\Omega'$.
\end{theorem}

\subsection{Uniqueness}\label{sec4.1}
We now recall where we are at. We have $h:\oD\to\oD$ a homeomorphic solution to (\ref{3.1}). In any local Jordan domain $\Omega\Subset\ID$, $\Omega'=h(\Omega)$ we have shown above in Theorem \ref{diffeoExistence} that there is a diffeomorphic solution $H:\Omega\to\Omega'$ to the {\em same} Alfors-Hopf equation and we may further assume that there $w_0\in\Omega$ and $a\in \partial \Omega$ with $H(w_0=h(w_0)$ and $H(a)=h(a)$. This is the usual two-point normalisation for uniqueness in elliptic problems.  However our equations are not elliptic. However our task in here is to prove that they are the same function. As noted earlier, we will deal with the zero and non-zero points of $\Phi$ in different ways. Recall $Z=\{w\in\ID:\Phi(w)=0\}$.

\subsubsection{Non-zero points of $\Phi$.}
We now consider a point $w_0\in\ID\setminus Z$. If we choose the neighbourhood $w_0\in\Omega\Subset\ID\setminus Z$ sufficiently small, then in $\Omega$ there is a conformal mapping $\psi$ which satisfies $\psi'=\sqrt{\Phi}$, where $\sqrt{\Phi}$ is any well-defined branch in $\Omega$. We consider $\tilde{h}=h\circ\psi^{-1}$. A calculation reveals this function satisfies
\begin{eqnarray}
\lefteqn{\exp(p\IK(w,\tilde{h}))\tilde{h}_w\overline{\tilde{h}_\wbar}\eta(\tilde{h}) }\nonumber \\ &= & [\exp(p\IK(w,h))h_w\overline{h_\wbar}\eta(h)](\psi^{-1})\cdot(\psi^{-1})'^2\; =\; 1.\label{4.1}
\end{eqnarray}
Furthermore, we may choose a disk $D$ in $\psi(\Omega)$ and consider $\tilde{h}=h\circ\psi^{-1}:D\to h\circ\psi^{-1}(D):=\tilde{\Omega}$. Next we consider the equation for $\tilde{f}=\tilde{h}^{-1}$ by substitution to see that
\begin{equation}\label{4.2}
\exp(p\IK(z,{\tilde{f}}))\frac{\overline{\tilde{f}_z\tilde{f}_\zbar}}{(J(z,\tilde{f}))^2}\eta=-1.
\end{equation}
In fact we only require  the right hand side to be a constant so as to make the equation homogeneous. By Theorem \ref{diffeoExistence}, there is a diffeomorphic solution $H:\overline{D}\to\overline{\tilde{\Omega}}$ to the same Ahlfors-Hopf equation (\ref{4.1}), and then its inverse $g=H^{-1}$ also satisfies (\ref{4.2}). We now have the following uniqueness result:\\

\begin{theorem}\label{Theo4}
Let $\ID$ be the unit disk and $\Omega$ be a Jordan domain. Let $f,g:\overline{\Omega}\to\oD$ both be finite distortion homeomorphisms which satisfy the equation
\begin{equation}\label{4.3}
\exp(p\IK(z,f))\frac{f_zf_\zbar}{J(z,f)^2}=\frac{1}{\eta},
\end{equation}
where $\eta\in C^\infty(\Omega)$, $\eta\geq1$ and there are points $z_0\in\Omega$, $w_0\in\ID$, $a\in\partial\Omega$, $b\in\IS$ such that $f(z_0)=g(z_0)=w_0$, $f(a)=g(a)=b$.  Then $f\equiv g$ in $\overline{\Omega}$.
\end{theorem}
Now we show (\ref{4.3}) can be regarded as an equation $f_\zbar=\mathcal{H}(z,f_z)$, and $g$ satisfies the same equation $g_\zbar=\mathcal{H}(z,g_z)$. We take the absolute values and consider the equation
\begin{equation}\label{4.4}
\exp\left(p\frac{x^2+y^2}{x^2-y^2}\right)\frac{xy}{(x^2-y^2)^2}=k,
\end{equation}
where $x,y,k$ are all non-negative, and $y\leq x$. We implicitly differentiate $y$ as a function of $x$ to obtain
\begin{equation}\label{4.5}
y'(x)=\frac{y}{x}\frac{3x^4+(4p-2)x^2y^2-y^4}{x^4+(4p+2)x^2y^2-3y^4}=t\frac{3+(4p-2)t^2-t^4}{1+(4p+2)t^2-3t^4},\quad t=\frac{y}{x}.
\end{equation}
\[
y''(x)=\frac{2 t \left(1-t^2\right)^2 \left[3 t^8-12 t^6+2(9+8p^2) t^4-12 t^2+3\right]}{x \left[-3 t^4+(4p+2) t^2+1\right]^3} > 0.
\]
Thus $y'(x)$ is increasing and in (\ref{4.5}) we also observe
\[
\frac{y}{x}\leq y'\leq3\frac{y}{x}.
\]
As the left-hand side of (\ref{4.4}) is strictly increasing for $y\in [0,x]$ and has range $[0,\infty]$ we have a well-defined function,  for fixed $k$,
\[ 
y = \mathcal{A}_k(x), \quad y<x. 
\]
This yields 
\begin{equation}\label{4.6}
 f_\zbar =\mathcal{A}_{\frac{1}{\eta(z)}} (|f_z|)\frac{\overline{f_z}}{|f_z|}:= \mathcal{H}(z,f_z).
\end{equation}
We seek the elliptic estimate and for this we need to compute
\begin{eqnarray*}
\left|\mathcal{A}_k(|\zeta|)-\mathcal{A}_k(|\xi|)\right| & = & \int_{|\xi|}^{|\zeta|} \mathcal{A}_k'(x) \; dx= \int_{|\xi|}^{|\zeta|}y'(x) dx .
\end{eqnarray*}
We see then that as $y'$ is increasing,
\begin{eqnarray}\label{4.7}
\left|\mathcal{A}_k(|\zeta|)-\mathcal{A}_k(|\xi|)\right|  & \leq &  \mathcal{A}_k'(|\zeta|)\big(|\zeta|-|\xi|),  \quad   \mbox{ if }|\zeta|\geq |\xi|.
\end{eqnarray}
With $v_k(x)=\mathcal{A}_k(x)/x$,
\begin{eqnarray}
\frac{|\mathcal{H}(z,\zeta)-\mathcal{H}(z,\xi)|}{|\zeta-\xi|}&=&\frac{\left|\mathcal{A}_k(|\zeta|)\frac{\zeta}{|\zeta|}-\mathcal{A}_k(|\xi|)\frac{\xi}{|\xi|}\right|}{|\zeta-\xi|}\notag\\
&=&\frac{|v_k(|\zeta|)\zeta-v_k(|\xi|)\xi|}{|\zeta-\xi|}.\label{4.8}
\end{eqnarray}

\begin{lemma}\label{maxineq}
Let $a,b\in\IR$, $a,b\geq0$, $z,w\in\IC$. Then
\[
\frac{|az-bw|}{|z-w|}\leq\max\left\{\frac{a|z|+b|w|}{|z|+|w|},\left|\frac{a|z|-b|w|}{|z|-|w|}\right|\right\}.
\]
\end{lemma}
\noindent{\bf Proof.} Let $|z|=t$, $|w|=s$, then for some angle $\theta\in[0,2\pi)$ we have
\[
\frac{|az-bw|^2}{|z-w|^2}=\frac{|at-bse^{i\theta}|^2}{|t-se^{i\theta}|^2}=\frac{a^2t^2+b^2s^2-2abst\cos\theta}{t^2+s^2-2st\cos\theta}:=F(\theta).
\]
We differentiate this with respect to $\theta$ to see
\[
\frac{d}{d\theta}F(\theta)=\frac{2  st(a -b )(bs^2 -at^2 )}{(t^2+s^2-2st\cos\theta)^2}\sin\theta.
\]
Now it entirely depends on the sign of the term $(a -b )(bs^2 -at^2 )$. If it is positive, then $F(\theta)$ gets maximum at $\theta=\pi$, which gives the first case; if it is negative, then $F(\theta)$ gets maximum at $\theta=0$, which gives the second case; if it is $0$, then $F(\theta)$ is a constant for $\theta\in[0,2\pi)$, and the equality holds.\hfill$\Box$

\bigskip

Now by (\ref{4.7}), (\ref{4.8}), and Lemma \ref{maxineq}, we have
\begin{eqnarray*}
\frac{|\mathcal{H}(z,\zeta)-\mathcal{H}(z,\xi)|}{|\zeta-\xi|}&\leq&\max\left\{\frac{\mathcal{A}_k(|\zeta|)+\mathcal{A}_k(|\xi|)}{|\zeta|+|\xi|},\frac{|\mathcal{A}_k(|\zeta|)-\mathcal{A}_k(|\xi|)|}{||\zeta|-|\xi||}\right\}\\
&\leq&\max\{v_k(|\zeta|),v_k(|\xi|),\mathcal{A}'_k(|\zeta|),\mathcal{A}'_k(|\xi|)\}.
\end{eqnarray*}
We now consider
\begin{equation}\label{4.9}
\mu:=\frac{|\mathcal{H}(z,f_z)-\mathcal{H}(z,g_z)|}{|f_z-g_z|}\leq\max\{v_\frac{1}{\eta}(|f_z|),v_\frac{1}{\eta}(|g_z|),\mathcal{A}'_\frac{1}{\eta}(|f_z|),\mathcal{A}'_\frac{1}{\eta}(|g_z|)\}.
\end{equation}
As $v_\frac{1}{\eta}(|f_z|)=|\mu_f|$ and $v_\frac{1}{\eta}(|g_z|)=|\mu_g|$, they are both in the exponential class. For the last two terms in (\ref{4.9}), recall our computation (\ref{4.5}). Then with $y(x)=\mathcal{A}_k(|\zeta|)$,
\[
y'(x)=t\frac{ 3  +(4p-2)t^2-t^4 }{1+(4p+2)t^2-3 t^4 },\quad t=\frac{y}{x}.  
\]
Then with $R_p(t) =\frac{1+y'^2}{1-y'^2}\Big/\frac{1+t^2}{1-t^2}$ we calculate
{\footnotesize \[
R_p(t) = \frac{1+t^2(13+8p+2(-7+4p(5+2p))t^2+2(-7+4p(-5+2p))t^4+(13-8 p)t^6+t^8)}{(1+t^2)(1+(-4+8p)t^2+2(3+8p^2) t^4-4(1+2p)t^6 +t^8)}.
\]}
We are only interested in the situation when $t\to1$ and so observe that $R_p(1)=1$ for any $p$.

\begin{center}
\scalebox{0.5}{\includegraphics*{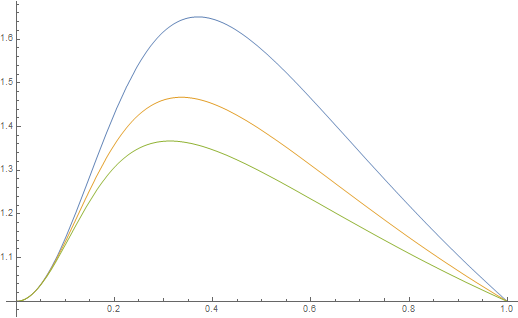}} \\
The graph of $R_p(t)$ for $p=2,3,4$. 
\end{center}
This gives 
\[ s\; \frac{1+y'^2}{1-y'^2}\leq p\; \frac{1+t^2}{1-t^2}+C(p,s) \]
 for any $0< s<p$.\\

Now in (\ref{4.9}), if we put $F=g-f$, then we find $F$ satisfies the Beltrami equation
\[
F_\zbar=\mu F_z,
\]
where for any $0<s<p$,
\[
\int_\Omega\exp\left(s\; \frac{1+|\mu|^2}{1-|\mu|^2}\right)<\infty.
\]
We thus have that $F$ is open if not a constant, see \cite{KKM2} or \cite[Theorem 20.4.9]{AIM}. Also the image of $\partial\Omega$ is the difference of two homeomorphisms from $\partial\Omega$ to $\IS$, which has \textit{total variation of the argument} less than $2\pi$, and this contradicts the openness of $F$. So we conclude that $F\equiv0$. See \cite[Section 9.2.2]{AIM} for the precise details here. This completes the proof of Theorem \ref{Theo4}.\\

We remark here that  the reason that we transformed the Ahlfors-Hopf differential $\Phi$ to $1$ to prove the uniqueness on the $f$ side is that is we tried this on the $h$ side, then there is a term $\eta(h)$; if $\Phi$ is not a constant, then it will become $\Phi(f)$ when we transform the equation to $f$ side. All the efforts we made were to remove the $f$ or $h$ term to make the equation homogeneous;  if we were consider the difference $g-f$ or $H-h$, we would have an equation in the form 
\[
F_\zbar=\mu F_z+\sigma F.
\]
Such equations requires higher global regularity of the functions to prove  uniqueness. In fact we will meet this situation now.

\subsubsection{Zero points of $\Phi$.}
We now come to the case $w_0\in Z$. As we have proved in the last section, $h$ is diffeomorphic in $\ID\setminus Z$, and since $Z$ is discrete, we can choose a neighbourhood $w_0\in\Omega\Subset\ID$ such that $h$ is diffeomorphic in $\overline{\Omega}\setminus\{w_0\}$. We may also assume the image $h(\Omega)=D$ is a disk. By Theorem \ref{diffeoExistence}, there is a $H:\overline{\Omega}\to D$, $H$ is diffeomorphic in $\Omega$, and $H$ satisfies the same Ahlfors-Hopf equation. Thus we can set the problem as follows:
\begin{theorem}\label{Theo5}
Let $\ID$ be the unit disk and $\Omega$ be a Jordan domain. Let $h,H:\overline{\Omega}\to\oD$ both be finite distortion homeomorphisms which satisfy the equation
\begin{equation}\label{4.10}
\Phi=\exp(p\IK(w,h))h_w\overline{h_\wbar}\eta(h),
\end{equation}
where $\eta\in C^\infty(\Omega)$, $\eta\geq1$ and there are points $w_0\in\Omega$, $z_0\in\ID$, $a\in\partial\Omega$, $b\in\IS$ such that $h(w_0)=H(w_0)=z_0$, $h(a)=H(a)=b$. Furthermore, $h$ is diffeomorphic in $\oD\setminus\{w_0\}$, $H$ is diffeomorphic in $\ID$. Then $h=H$ in $\overline{\Omega}$.
\end{theorem}

\bigskip

Here we may write (\ref{4.10}) as an equation
\[
h_\wbar=\mathcal{B}(w, h, h_w).
\]
We separate the domain $\Omega$ into $\Omega_r$ and $\Omega\setminus\Omega_r$, where $w_0\in\Omega_r\Subset\Omega$ and $r$ will be small. We first consider the points $w\in\Omega_r$, since $H$ is diffeomorphic in $\overline{\Omega_r}$, we have $|H_w|\geq\epsilon>0$. Thus
\begin{eqnarray*}
|h_\wbar-H_\wbar|&=&\left|\mathcal{B}(w,h,h_w)-\mathcal{B}(w,H,H_w)\right|\\
&\leq&\left|\mathcal{B}(w,h,h_w)-\mathcal{B}(w,h,H_w)\right|+\left|\mathcal{B}(w,h,H_w)-\mathcal{B}(w,g,H_w)\right|\\
&=&|\mu||h_w-H_w|+|\sigma||h-H|.
\end{eqnarray*}
where
\[
\mu=\frac{\mathcal{B}(w,h,h_w)-\mathcal{B}(w,h,H_w)}{h_w-H_w},\quad\sigma=\frac{\mathcal{B}(w,h,H_w)-\mathcal{B}(w,H,H_w)}{h-H}.
\]
We need estimates on both $\mu$ and $\sigma$. Note that for $\mu$ the terms $\mathcal{B}(w,h,h_w)$ and $\mathcal{B}(w,h,H_w)$ have the same second variable $h$, so (\ref{3.17}) in Lemma \ref{qrlemma} applies, and so yields
\[
|\mu|\leq\frac{1}{2}(1+v(|H_w|))<1,
\]
where $v(|H_w|)=\frac{|\mathcal{B}(w,h,H_w)|}{|H_w|}<1$. For the estimate on $\sigma$, we consider the equation for $\mathcal{A}=\mathcal{A}_k(x)$:
\[
\exp\left(p\frac{x^2+\mathcal{A}^2}{x^2-\mathcal{A}^2}\right)\mathcal{A}x=k.
\]
Put $\log$ on both sides and differentiate  with respect to  $k$ to obtain
\[
\frac{\partial \mathcal{A}}{\partial k}\left(\frac{4p\mathcal{A}x^2}{(x^2-\mathcal{A}^2)^2}+\frac{1}{\mathcal{A}}\right)=\frac{1}{k}.
\]
Here
\[
\frac{4p\mathcal{A}x^2}{(x^2-\mathcal{A}^2)^2}+\frac{1}{\mathcal{A}}\geq\frac{1}{\mathcal{A}}.
\]
So
\[
\frac{\partial\mathcal{A}}{\partial k}\leq\frac{\mathcal{A}}{k}=\frac{1}{\exp\left(p\frac{x^2+\mathcal{A}^2}{x^2-\mathcal{A}^2}\right)x}\leq\frac{1}{x}.
\]
Note the last estimate does not depend on $k$. Also from (\ref{4.10}), $\mathcal{B}(w,h,H_w)$ and $\mathcal{B}(w,H,H_w)$ have the same argument. Thus for some $k_1$ between $\frac{|\Phi|}{\eta(h)}$ and $\frac{|\Phi|}{\eta(H)}$,
\begin{eqnarray*}
\left|\mathcal{B}(w,h,H_w)-\mathcal{B}(w,H,H_w)\right|&=&\left|\mathcal{A}_{\frac{|\Phi|}{\eta(h)}}(|H_w|)-\mathcal{A}_{\frac{|\Phi|}{\eta(H)}}(|H_w|)\right|\\
&=&\frac{\partial\mathcal{A}_k(|H_w|)}{\partial k}\Big|_{k=k_1}\left|\frac{|\Phi|}{\eta(h)}-\frac{|\Phi|}{\eta(H)}\right|\notag\\
&\leq&\frac{C}{|H_w|}|h-H|\leq\frac{C}{\epsilon}|h-H|.
\end{eqnarray*}
This proves $\sigma\in L^\infty(\Omega_r)$. In $\Omega\setminus\Omega_r$, we have $|h_w|\geq\epsilon>0$. Thus we can interchange the roles of $h$ and $H$ and get
\[
|h_\wbar-H_\wbar|\leq|\mu||h_w-H_w|+|\sigma||h-H|,
\]
where
\[
\mu=\frac{\mathcal{B}(w,H,H_w)-\mathcal{B}(w,H,h_w)}{H_w-h_w},\quad\sigma=\frac{\mathcal{B}(w,H,h_w)-\mathcal{B}(w,h,h_w)}{H-h},
\]
and the same argument follows. We conclude that for the difference $F=H-h$,
\[
F_\wbar=\mu F_w+\sigma F,
\]
where $|\mu|\leq k<1$ in $\Omega$, and $\sigma\in L^\infty(\Omega)$. This is now enough to prove that $F\equiv0$ in $\Omega$ following the argument of  \cite[Section 9.2.2]{AIM}.

\section{The limiting regimes, $p\to\infty$ and $p\to 0$.}
In this section we consider what happens as $p$ varies in the $\exp(p)$ problems. In particular, we wish to recover the harmonic mappings as $p\to0$ and the Teichm\"uller mappings as $p\to\infty$. Notice that for each $0<p<\infty$, we have already shown that there exists a unique diffeomorphic minimiser $f_p$, which are thus locally quasiconformal - and quasiconformal on any relatively compact domain which we will assume $R$ is.\\

To consider what happens as $p\to0$ and $p\to\infty$, we need to redefine our functionals as
\[
\mE_p(f)=\frac{1}{p}\log\frac{1}{|R|}\int_R\exp(p\IK(z,f))\;d\sigma_R,
\]
where $|R|$ denotes the area of $R$. In this way we can extend the definitions to $p=0$ and $p=\infty$ as
\[
\mE_0(f)=\frac{1}{|R|}\int_R\IK(z,f)\;d\sigma_R,
\]
\[
\mE_\infty(f)=\sup_{z\in R}\IK(z,f).
\]
In this section we will prove the following:
\begin{theorem}
Let $f_p$ be the minimiser of $\mE_p(f)$, $h_p=f_p^{-1}$. Let $0\leq p_0\leq\infty$, $p\to p_0$. Then
\begin{itemize}
\item The function $p\mapsto\mE_p(f_p)$ is a continuous non-decreasing function for $p\in[0,\infty]$.
\item If $0<p_0<\infty$, then both $f_p\to f_{p_0}$, $h_p\to h_{p_0}$ uniformly with all derivatives.
\item In the case $p_0=0$ we have $h_p\to h_0$ uniformly, where $h_0$ is a harmonic mapping; 
\item In the case $p_0=\infty$, then both $f_p\to f_\infty$ and $h_p\to h_\infty$ uniformly, where $h_\infty$ is a Teichm\"uller mapping.
\end{itemize}
\end{theorem}

\subsubsection{$0\leq p_0<\infty$ case.}
Let $p\to p_0$, $0\leq p_0<\infty$. Let $f_p$ be the unique diffeomorphic minimiser of $\mE_p$, $f$ be the weak limit of $f_p$ in $W^{1,2}(R)$, $f_{p_0}$ be the unique diffeomorphic minimiser of $\mE_{p_0}$. Since $f_{p_0}$ is quasiconformal, $\mE_p(f_{p_0})$ is finite for all $p$, and $\mE_p(f_{p_0})\to\mE_{p_0}(f_{p_0})$ for $p\to p_0$. Then,
\[
\mE_{p_0}(f)\leq\lim_{p\to p_0}\mE_p(f_p)\leq\lim_{p\to p_0}\mE_p(f_{p_0})=\mE_{p_0}(f_0),
\]
where the first inequality follows from the polyconvexity, and the second is because each $f_p$ is the minimiser of $\mE_p$. However, since $f_0$ is the unique minimiser of $\mE_{p_0}$, we then have $f=f_0$ and
\begin{equation}\label{5.1}
\lim_{p\to p_0}\mE_p(f_p)=\mE_{p_0}(f_0).
\end{equation}
Note this works for any $p_0\in[0,\infty)$. In fact, for $0<p_0<\infty$, we have proved the minimiser $f_{p_0}$ is a quasiconformal diffeomorphic mapping; for $p_0=0$, we also have a quasiconformal diffeomorphic minimiser, since its inverse is harmonic, see \cite{SY}.\\

There are several ways to see the stronger convergence for $0<p_0<\infty$. Firstly, it follows from (\ref{5.1}) and  the Radon-Riesz property that $f_p\to f_{p_0}$ strongly in $W^{1,2}(R)$, $h_p\to h_{p_0}$ strongly in $W^{1,2}(S)$, see \cite{MY3}. \\

For the convergence of the higher order derivatives, recall in \S 2.3, using the inner variational equations, we found a sequence of functions $F^p$, which satisfy
\[
F^p_z=\exp(p\IK(z,f))\frac{2p\overline{\mu_f}}{1-|\mu_f|^2}\eta,\quad F^p_\zbar=(\exp(p\IK(z,f))-e^p)\eta+\sigma^p.
\]
These satisfy the elliptic equations
\[
F^p_\zbar=\mathcal{A}_p\left(\frac{|F^p_z|}{\eta}\right)\eta+\sigma^p,
\]
where $\mathcal{A}_p'\leq k_p<1$, which is uniform with $p$. Following our arguments as in \S 2.3, this sequence also converges with all derivatives as $p\to p_0$, for $0<p_0<\infty$. Thus $\mu_p\to\mu$, and then $f_p\to f$, $h_p\to h$, all with arbitrary order of derivatives.\\

Another way is through the $\mu$ equation (\ref{2.14}). We have seen that each $\mu_p$ satisfies the equations
\[
(\mu_p)_z=\nu_p(\mu_p)_\zbar+\phi_p,
\]
where $\phi_p$ are uniformly bounded in $L^\infty(R)$, and
\[
\frac{1}{p}\log\frac{1}{|R|}\int_R\exp\left(c(p)\sqrt{\frac{1+|\nu_p|^2}{1-|\nu_p|^2}}\right) \;d\sigma_R \leq\frac{1}{p}\log\frac{1}{|R|}\int_R\exp(p\IK(z,f_p)) \; d\sigma_R
\]
are uniformly bounded, where $c(p)$ is a constant depending on $p$. This implies that $\mu_p$ are uniformly bounded in $W^{1,s}(S)$ for some $1<s<2$ (cf. \cite[Theorem 20.4.3]{AIM}), and then it follows from the good approximation lemma that $\mu_p$ converge to $\mu_f$, see \cite[Lemma 5.3.5]{AIM}.\\

Note neither of these methods works for the case $p\to\infty$. In fact no equation would be uniformly elliptic in the general case, because the limit is a Teichm\"uller mapping, which has constant $|\mu|$. Thus
\[
0=(|\mu|^2)_z=\mu_z\mubar+\mu\overline{\mu_\zbar}.
\]
To solve this problem, recall in \S \ref{sec4.1} we introduced a way to turn the Beltrami coefficient into a real function. We will exploit this in the next section.

\subsubsection{$p_0=\infty$ case.} 
As we stated in \S 1.1 and 1.2, Ahlfors has already given a proof for Teichm\"uller's theorem using the minimisers of $L^p$ mean distortion in a homotopy class $[f_0]$ and letting  $p\to\infty$, 
\[ \min \;\Big\{ \int_\Sigma \IK(z,f)^p d\sigma : f\in [f_0] \Big\} \]
\cite{Ah3}. Here we give a more computational and  easier proof since  we now have the diffeomorphic $\exp(p)$ minimisers, which have far more regularity than the $L^p$ minimisers that Ahlfors had. We also have decades of new technology to work with.\\

We first note that $\mE_{p_0}(f_p)\leq\mE_p(f_p)$ are uniformly bounded for all $p\geq p_0>1$, so $f_p$ are uniformly bounded in $W^{1,2}(R)$; $h_p$ are all uniformly bounded in $W^{1,2}(S)$ as $p\to\infty$, so there exists a homeomorphism $f_\infty$ and $h_\infty=f_\infty^{-1}$ such that $f_p\to f_\infty$ and $h_p\to h_\infty$ uniformly. Now using polyconvexity, we have
\[
\mE_\infty(f_\infty)\leq\lim_{p\to\infty}\mE_p(f_\infty)\leq\lim_{p\to\infty}\mE_p(f_p).
\]
If the inequality holds, then for sufficiently large $p$ we have $\mE_p(f_\infty)<\mE_p(f_p)$, which contradicts that $f_p$ is the minimiser of $\mE_p$. This implies that $f_\infty$ must be a minimiser of $\mE_\infty$, and we have
\[
\mE_\infty(f_\infty)=\lim_{p\to p_0}\mE_p(f_p).
\]

We lift the functions to the planar disk. Then the holomorphic Ahlfors-Hopf differential equation reads as
\[
\Phi_p=c(p)\exp(p\IK(w,h))h_w\overline{h_\wbar}\eta(h),
\]
where $\eta(z)=\frac{1}{(1-|z|^2)^2}$, $c(p)$ is a constant chosen so that $\|\Phi_p\|_{L^1(\tilde{\mP})}=1$, for any fundamental polyhedron $\tilde{\mP}$ of $S$. Following our earlier argument, by the Riemann-Roch theorem, we have $\Phi_p\to\Phi\not\equiv0$ uniformly, where $\Phi$ is also a holomorphic mapping.\\

Now let $Z$ be the zero set of $\Phi$, and $w_0\in\ID\setminus Z$. As $\Phi_p\to\Phi$ uniformly, we may choose a small neighbourhood $\Omega$ of $w_0$ such that for sufficiently large $p$, there is a conformal mapping $\Psi_p$ that satisfies $\Psi_p'=\sqrt{\Phi_p}$. Then $\tilde{h}_p=h_p\circ\Psi_p^{-1}$ satisfies
\begin{equation}\label{5.2}
c(p)\exp(p\IK(w,\tilde{h}_p))(\tilde{h}_p)_w\overline{(\tilde{h}_p)_\wbar}\eta(\tilde{h})=1.
\end{equation}
We drop the sup and subscripts. Note this means $\mu=\mu_{h_p}$ is real positive. We choose any $\Omega'\Subset f_\infty(\Omega)$, and then transform this to the ``$f$ side'' to get
\[
\nu:=\mu(f)=-\frac{f_\zbar}{\overline{f_z}}.
\]
Now $\nu(z)$ is also a real positive function, which satisfies
\begin{equation}\label{5.3}
\nu_z=-\left(\frac{f_\zbar}{\overline{f_z}}\right)_z=-\frac{1}{\overline{f_z}^2}(f_{z\zbar}\overline{f_z}-f_\zbar\overline{f_{z\zbar}})=-\frac{1}{\overline{f_z}}(f_{z\zbar}+\nu\overline{f_{z\zbar}}).
\end{equation}
Then (\ref{5.2}) reads as
\[
c(p)\exp\left(p\frac{1+\nu^2}{1-\nu^2}\right)\frac{\nu}{(1-\nu^2)^2}\eta=|f_z|^2.
\]
Now we can compute
\begin{eqnarray*}
\log c(p)+p\frac{1+\nu^2}{1-\nu^2}+\log\nu-2\log(1-\nu^2)+\log\eta&=&\log|f_z|^2,\\
\left(\frac{4p\nu}{(1-\nu^2)^2}+\frac{1}{\nu}+\frac{4\nu}{1-\nu^2}\right)\nu_z+\frac{\eta_z}{\eta}&=&\frac{1}{|f_z|^2}(f_{zz}\overline{f_z}+f_z\overline{f_{z\zbar}}),\\
-\frac{1+(4p+2)\nu^2-3\nu^4}{\nu(1-\nu^2)^2}(f_{z\zbar}+\nu\overline{f_{z\zbar}})+\frac{\eta_z}{\eta}\overline{f_z}&=&\xi f_{zz}+\overline{f_{z\zbar}},
\end{eqnarray*}
where $\xi=\frac{\overline{f_z}}{f_z}$. We write
\[
\alpha=1+(4p+2)\nu^2-3\nu^4,\quad\beta=\nu(1-\nu^2)^2,\quad\phi=\frac{\eta_z}{\eta}\overline{f_z}.
\]
Then the last equation reads as
\[
\alpha f_{z\zbar}+(\nu\alpha+\beta)\overline{f_{z\zbar}}+\beta\xi f_{zz}=\beta\phi.
\]
We solve for $f_{z\zbar}$ in the last equation, and obtain
\begin{equation}\label{5.4}
\Gamma f_{z\zbar}=A\xibar\overline{f_{zz}}-B\xi f_{zz}+\psi,
\end{equation}
where
\[
\Gamma=\alpha^2-(\nu\alpha+\beta)^2,\quad A=(\nu\alpha+\beta)\beta,\quad B=\alpha\beta,
\]
\[
\psi=\alpha\beta\phi-(\nu\alpha+\beta)\beta\overline{\phi}.
\]
Note then as $p\to\infty$,
\begin{eqnarray*}
\frac{A+B}{\Gamma}&=&\frac{\beta}{(1-\nu)\alpha-\beta}=\frac{\nu(1+\nu)(1-\nu^2)}{1+(4p+2)\nu^2-3\nu^4-\nu(1+\nu)(1-\nu^2)}\\
&=&\frac{\nu(1+\nu)(1-\nu^2)}{4p\nu^2+(1-\nu+2\nu^2)(1-\nu^2)}\leq \frac{c}{\sqrt{p}},\\
\frac{|\psi|}{\Gamma}&\leq&\frac{\beta}{(1-\nu)\alpha-\beta}|\phi|\leq\frac{c}{\sqrt{p}}\frac{|\eta_z|}{\eta}|f_z|\in L^2(\Omega').
\end{eqnarray*}
This makes (\ref{5.4}) an elliptic equation and thus $f_p$ are uniformly bounded in $W^{2,2}(\Omega')$ as $p\to\infty$. In particular, $Df_p\to Df_\infty$ pointwise, so $\nu_{f_p}\to\nu_{f_\infty}$ pointwise, and then $\nu_{f_\infty}=\mu_{h_\infty}(f_\infty)$ by the good approximation lemma. Now let $p\to\infty$, then (\ref{5.4}) reads as
\[
(f_\infty)_{z\zbar}=0.
\]
Then by (\ref{5.3}),
\[
(\nu_{f_\infty})_z=0.
\]
Thus $\mu_{h_\infty}=\nu_{f_\infty}(h_\infty)$ is a constant.\\

We now put back the sup and subscripts. For the sequence $\tilde{h}_p$, we have proved the limit function $\tilde{h}_\infty$ has $\mu_{\tilde{h}_\infty}=k$, a constant. Since $\Phi_p\to\Phi$, we also have $\Psi_p\to\Psi$, where $\Psi$ is also a non-zero holomorphic function and $\Psi'=\sqrt{\Phi}$. Thus $h_p\to h_\infty=\tilde{h}_\infty\circ\Psi$, and
\[
\mu_{h_\infty}=\frac{(h_\infty)_\wbar}{(h_\infty)_w}=\frac{(\tilde{h}_\infty)_\wbar(\Psi)\overline{\Psi'}}{(\tilde{h}_\infty)_w(\Psi)\Psi'}=k\frac{\sqrt{\overline{\Phi}}}{\sqrt{\Phi}}=k\frac{\overline{\Phi}}{|\Phi|}.
\]
This equation holds at every point $w\in S\setminus Z$, where $Z$ is a discrete set, since $\Phi$ is holomorphic, we have recovered Teichm\"uller's equation (\ref{teich1}).

\subsection{An example.} In the argument above we have $h_p \to h_\infty$,  $h_p$ are all diffeomorphisms and $h_\infty$ is a Teichm\"uller mapping.  The convergence is strong in $W^{1,2}(\ID)$,  and actually this convergence is with all derivatives away from the zero set of $\Phi$. Notice that $h_\infty$ cannot be a diffeomorphism at a zero of $\Phi$ as $\mu_{h_\infty}$ will be discontinuous there. This outcome is not surprising.  Let $0< k <1$ and consider the Beltrami equation
\[ g_\zbar = k |z|^{1/p} \frac{\bar{z}}{z} \; g_z \]
and solutions $g_p:\ID \to \ID$, $g_p(0)=0$. The Beltrami coeffficients $\mu_{g_p}=k |z|^{1/p} \frac{\bar{z}}{z}$ are continuous, each $g_p$ is a k-quasiconformal diffeomorphism and $\mu_{g_p} \to \frac{k \bar{z}}{z}$ with $p\to\infty$ locally uniformly in $\ID\setminus\{0\}$  so the good approximation theorem gives $g_\infty$ quasiconformal and a Teichm\"uller mapping,  which is not a diffeomorphism near $0$. 

This shows we cannot expect locally uniform $W^{2,2}$ bounds in our problem near $\{\Phi=0\}$.
\section{General boundary value problems}
In this section we give some brief remarks on the general boundary value problems for functions of exponentially finite distortion. Consider the energy functional
\[
E_p(f):=\int_\ID\exp(p\IK(z,f))\;\eta(z)dz,
\]
where $\eta$ is any positive weight. Now the problem becomes: Given a homeomorphism $f_0:\oD\to\oD$ with $E_p(f_0)<\infty$, we wish to find a minimiser of $E_p$ in the space
\[
{\cal F}:=\left\{f\in W_{loc}^{1,1}(\ID):E_p(f)<\infty,f|_\IS=f_0|_\IS, f:\oD\to\oD\mbox{ is a homeomorphism}\right\}.
\]
In fact, almost everything we have done in the surface case also works for this problem, apart from the existence of Ahlfors-Hopf differential, where the Riemann-Roch theorem was applied. In fact, if a minimiser $f$ is variational, then $f$ and its inverse $h$ satisfy the following equations
\begin{equation}\label{6.1}
\int_\ID\exp(p\IK(z,f))(\eta\varphi)_z \;dz=\int_\ID\exp(p\IK(z,f))\frac{2p\overline{\mu_f}}{1-|\mu_f|^2}\eta\varphi_\zbar\;dz,\quad\forall\varphi\in C_0^\infty(\ID);
\end{equation}
\begin{equation}\label{6.2}
\int_\ID\exp(p\IK(w,h))h_w\overline{h_\wbar}\eta(h)\psi_\wbar\;dw=0,\quad\forall\psi\in C_0^\infty(\ID).
\end{equation}
This variationality is not automatically true, \cite{MY2}. However, if the Ahlfors-Hopf equation holds, then the same methods we have already developed work and we have $f$ is diffeomorphic.\\

It is a simple observation that  if $\exp(q\IK(z,f))\eta\in L^1_{loc}(\ID)$ for some $q>p$, then $f$ is variational. Alternatively, to get the holomorphic Ahlfors-Hopf differential
\begin{equation}\label{6.3}
\Phi=\exp(p\IK(w,h))h_w\overline{h_\wbar}\eta(h),
\end{equation}
we can consider the truncated problem (similar to what we did in \S 2.2), and get the sequence
\begin{equation}\label{6.4}
\Phi_N=\sum_{n=0}^{N-1}\frac{p^n\IK^n(w,h_N)}{n!}(h_N)_w\overline{(h_N)_\wbar}\eta(h_N).
\end{equation}
Thus the problem of the existence of an Ahlfors-Hopf differential is implied by the convergence of this sequence to something nonzero. We record this as follows:
\begin{theorem}
In the space ${\cal F}$, there exists a minimiser of $E_p$. Furthermore, $f$ is diffeomorphic of $\ID$, if either
\begin{enumerate}
\item There exists a $q>p$ such that $\exp(q\IK(z,f))\eta\in L^1_{loc}(\ID)$;
\item The sequence $\Phi_N$ in (\ref{6.4}) is a Hamilton sequence, i.e.
\[
\frac{\Phi_N}{\|\Phi_N\|_{L^1(\ID)}}\not\to0.
\]
\end{enumerate}
\end{theorem}
The reader may see that the same results as described earlier work in the limiting regimes $p=0$ and $p=\infty$.  Of course for $p=0$ one may simply take the Poisson extension of the boundary values.

\medskip

GJM Institute for Advanced Study, Massey University, Auckland, New Zealand \\
email: G.J.Martin@Massey.ac.nz

CY Research Center for Mathematics and Interdisciplinary Sciences, Shandong University, 266237, Qingdao and Frontiers Science Center for Nonlinear Expectations, Ministry of Education, P. R. China \\
email: c.yao@sdu.edu.cn
\end{document}